\documentclass[a4paper]{amsart}
\usepackage{amssymb,a4wide,amsthm,amsmath,amsfonts}
\usepackage{mathrsfs}
\usepackage[usenames,dvipsnames]{xcolor}
\usepackage{tikz}
\usepackage{enumitem}
\usepackage{float}
\usepackage{mathtools}
\usepackage{thmtools}
\usepackage{bm}

\theoremstyle{plain}
\declaretheorem[name=Proposition, numberwithin=section]{proposition}
\declaretheorem[name=Theorem, sibling=proposition]{theorem}
\declaretheorem[name=Lemma, sibling=proposition]{lemma}
\declaretheorem[name=Corollary, sibling=proposition]{corollary}
\declaretheorem[name=Example,sibling=proposition]{example}
\declaretheorem[name=Counterexample,sibling=proposition]{counterexample}

\theoremstyle{definition}
\declaretheorem[name=Remark, sibling=proposition]{remark}
\declaretheorem[name=Notation, sibling=proposition]{notation}
\declaretheorem[name=Definition, sibling=proposition]{definition}
\declaretheorem[name=Convention,sibling=proposition]{convention}
\declaretheorem[name=Assumption,sibling=proposition]{assumption}

\usepackage[colorlinks, linkcolor=blue]{hyperref}

\newcommand{\vertiii}[1]{{\left\vert\kern-0.25ex\left\vert\kern-0.25ex\left\vert #1 \right\vert\kern-0.25ex\right\vert\kern-0.25ex\right\vert}}
\newcommand\numberthis{\addtocounter{equation}{1}\tag{\theequation}}
\renewcommand{\d}{\mathrm{d}}
\DeclareMathOperator{\E}{\mathbb{E}}
\DeclareMathOperator{\leom}{\lesssim_{\,\omega}}
\DeclareMathOperator{\geom}{\gtrsim_{\,\omega}}

\makeatletter
\@namedef{subjclassname@2020}{\textup{2020} Mathematics Subject Classification}
\makeatother

\begin{document}
\title[] {Besov-Orlicz Path regularity of non-Gaussian processes}
\author{Petr \v{C}oupek}
\address{Charles University\\
Faculty of Mathematics and Physics\\
Sokolovsk\'{a}~83\\
186~75\\
Prague~8\\
Czech Republic}
\email{coupek@karlin.mff.cuni.cz}

\author{Martin Ondrej\'{a}t}
\address{Czech Academy of Sciences\\
Institute of Information Theory and Automation\\
Pod Vod\'{a}renskou v\v{e}\v{z}\'{i}~4\\
180~00\\
Prague~8\\
Czech Republic}
\email{ondrejat@utia.cas.cz}



\begin{abstract}
In the article, Besov-Orlicz regularity of sample paths of stochastic processes that are represented by multiple integrals of order $n\in\mathbb{N}$ is treated. We give sufficient conditions for the considered processes to have paths in the exponential Besov-Orlicz space \[B_{\varPhi_{2/n},\infty}^\alpha(0,T)\qquad \mbox{with}\qquad \varPhi_{2/n}(x)=\mathrm{e}^{x^{2/n}}-1.\]
These results provide an extension of what is known for scalar Gaussian stochastic processes to stochastic processes in an arbitrary finite Wiener chaos. As an application, the Besov-Orlicz path regularity of fractionally filtered Hermite processes is studied. But while the main focus is on the non-Gaussian case, some new path properties are obtained even for fractional Brownian motions.
\end{abstract}

\maketitle

\section{Introduction}

It is well-known that the paths of the Wiener process belong to the Besov-Orlicz space $B_{\varPhi_2,\infty}^{1/2}(0,T)$ where $\varPhi_2(x)=\mathrm{e}^{x^2}-1$. The original proof of this result in \cite{Cie93} relies on intricate equivalences for Besov norms but a different proof is also available in \cite{HytVer08}.
From this result, one immediately obtains, for example, that Brownian paths belong to both the Besov space $B_{p,\infty}^{1/2}(0,T)$ for all $p\in [1,\infty)$ and the modulus H\"older space $C^{|r\log r|^{1/2}}([0,T])$ although historically, these two results came first; see \cite{Cie91} and \cite{Lev37}.

\bigskip

There is a number of generalizations of this result in various directions. In \cite{OndSimKup18}, it is shown that any continuous local martingale with Lipschitz continuous quadratic variation as well as solutions to stochastic differential equations with locally bounded non-linearities have paths in the space $B_{\varPhi_2,\infty}^{1/2}(0,T)$. In \cite{OndVer20}, it is shown that such regularity is also retained by stochastic convolutions (with values in $2$-smooth Banach spaces) and the result is also shown for strong solutions to stochastic $p$-Laplace systems in \cite{Wich21}. In \cite{Ver09}, it is shown that the fractional Brownian motion of Hurst parameter $\alpha\in (0,1)$ has paths in the Besov-Orlicz space $B_{\varPhi_2,\infty}^{\alpha}(0,T)$. There are also results on bifractional Brownian motions \cite{BouNach20} and L\'evy processes \cite{AziFaUns18, FagFalUns17, FagUnsWar17}.

\bigskip

The purpose of this article is to extend the results of \cite{Ver09} to non-Gaussian stochastic processes. In particular, we consider stochastic processes that fit the following scheme: Let $H$ be a real separable Hilbert space and $W=\{W(h)\}_{h\in H}$ an $H$-isonormal Gaussian process. For $n\in\mathbb{N}_0$, denote the multiple divergence operator of order $n$ by $W_n$. Let $G=\{G(t)\}_{t\in [0,T]}$, $T>0$, be the (jointly measurable) process represented by $G(t)=W_n(A_t)$ with $A_t\in H^{\otimes_2^n}$, $t\in [0,T]$. For example, one can consider the family of fractional Brownian motions, Hermite (and, in particular, Rosenblatt) processes, or the so-called fractionally filtered generalized Hermite processes; see \cite{BaiTaq14}.

\subsection{Main results}

A main result of the present article is \autoref{main_thm} which provides sufficient conditions under which process $G$ has paths in the Besov-Orlicz space \[B^{\alpha}_{\varPhi_{2/n},\infty}(0,T)\quad\mbox{ where }\quad \varPhi_{2/n}(x)=\mathrm{e}^{x^{2/n}}-1.\] These conditions are formulated in terms of the kernel $A_t$. In the scalar case, \autoref{main_thm} extends \cite[Theorem 5.1]{Ver09}, where Gaussian processes are considered, to a non-Gaussian setting. We thus make precise how the order of the Wiener chaos in which process $G$ lives influences the regularity of its paths. In particular, it is clear that the higher the order of the Wiener chaos, the worse regularity of paths we get. Moreover, in \autoref{thm:main_thm_2}, the result is refined and the precise pathwise behavior of the integral of the increments is obtained. Our results cover Gaussian processes (e.g. standard and fractional Brownian motions) but also non-Gaussian processes (e.g. higher-order Hermite processes). For example, we show that the Rosenblatt process with Hurst parameter $\alpha\in (1/2,1)$ has paths in the space $B_{\varPhi_{1},\infty}^\alpha(0,T)$. 

\subsection{Proof method and invalidity of Gebelein's inequality in higher-order Wiener chaoses} 

In order to establish such Besov-Orlicz regularity results, one would hope to proceed as in \cite{Ver09} (or \cite{BouNach20}). The proofs there rely on Gebelein's inequality \cite{Geb41} (see also \cite{BesCie06}):

\begin{theorem}[Gebelein's inequality]
Let $(\xi,\eta)^\top$ be a centered Gaussian vector in $\mathbb{R}^2$ with $\mathbb{E}\xi^2=\mathbb{E}\eta^2=1$. Then the inequality 
	\begin{equation}
	\label{eq:Gebelein}
	 |\mathbb{E} f(\xi)g(\eta)| \leq |\rho_{\xi,\eta}| \left[\mathbb{E} f(\xi)^2\right]^\frac{1}{2}\left[\mathbb{E} g(\eta)^2\right]^\frac{1}{2},
	 \end{equation}
where $\rho_{\xi,\eta}$ is the correlation coefficient between $\xi$ and $\eta$, holds for all functions $f,g:\mathbb{R}\rightarrow\mathbb{R}$ such that $\mathbb{E} f(\xi)^2<\infty$, $\mathbb{E} g(\eta)^2<\infty$, and $\mathbb{E}f(\xi)=\mathbb{E}g(\eta)=0$. 
\end{theorem}

However, while Gebelein's inequality can be used in a Gaussian setting, in our case, the process $G(t)=W_n(A_t)$, $t\in [0,T]$, is not Gaussian for $n\geq 2$. It is therefore natural to ask whether some analogue of Gebelein's inequality, that could be used for investigating Besov-Orlicz regularity of the paths of $G$, holds even in higher-order Wiener chaoses. There is the following partial result.

\begin{example}
Let $W$ be an $L^2(0,\infty)$-isonormal Gaussian process defined on $(\Omega,\mathcal{F},\mathbb{P})$ where the $\sigma$-algebra $\mathcal{F}$ is generated by $W$. Let $\xi=W_2(h_{\xi})$ for some symmetric function $h_{\xi}\in L^2((0,\infty)^2)$ with $\|h_{\xi}\|_{L^2((0,\infty)^2)}^2=1/2$, and let $\eta=W_2(h_{\eta}\otimes h_{\eta})$ for some function $h_\eta\in L^2(0,\infty)$ with $\|h_\eta\|_{L^2(0,\infty)}^4=1/2$. Then both $\xi$ and $\eta$ belong to the second Wiener chaos of $W$, $\mathbb{E}\xi=\mathbb{E}\eta=0$, and $\mathbb{E}\xi^2=\mathbb{E}\eta^2=1$. Moreover, the inequality
	\begin{equation}
	\label{eq:Gebelein_2}
	 |\mathbb{E}\,\xi g(\eta)| \leq |\rho_{\xi,\eta}|\, [\mathbb{E}\,\xi^2]^{\frac{1}{2}}[\mathbb{E}g(\eta)^2]^{\frac{1}{2}},
	 \end{equation}
where $\rho_{\xi,\eta}$ is the correlation coefficient between $\xi$ and $\eta$, holds for any function $g:\mathbb{R}\rightarrow\mathbb{R}$ such that $\mathbb{E}g(\eta)^2<\infty$.
\end{example}

\begin{proof}[Proof of inequality \eqref{eq:Gebelein_2}.]
By \cite[Proposition 3]{NuaUstZak88}, we have that $\mathbb{E}[\xi|\eta] = K\eta$ holds almost surely with 
	\[K=\frac{\langle h_\xi,h_{\eta}\otimes h_{\eta}\rangle_{L^2((0,\infty)^2)}}{\|h_\eta\otimes h_{\eta}\|_{L^2((0,\infty)^2)}^2}=\frac{\mathbb{E}[\xi\eta]}{\mathbb{E}\eta^2}.\]
Then we have 
	\begin{equation*}
		|\mathbb{E}\xi g(\eta)|  = |\mathbb{E}\left(\mathbb{E}\,[\xi|\eta]\,g(\eta)\right)| = |\mathbb{E} [K\eta g(\eta)]| \leq |K|\,[\mathbb{E} \eta^2]^\frac{1}{2} [\mathbb{E}g(\eta)^2]^\frac{1}{2}
	\end{equation*}
by using the Cauchy-Schwarz inequality and the claim follows.
\end{proof}

On the other hand, there is also the following negative result. 

\begin{counterexample}
Let $W$ be an $\mathbb{R}^2$-isonormal Gaussian process defined on $(\Omega,\mathcal{F},\mathbb{P})$ where the $\sigma$-algebra $\mathcal{F}$ is generated by $W$. Let $\{e_1,e_2\}$ be an orthonormal basis of $\mathbb{R}^2$ and set $X_i=W(e_i)$, $i=1,2$. Then both $X_1$ and $X_2$ are standard Gaussian random variables. Define $\xi=\frac{1}{\sqrt{2}}(X_1^2-1)$ and $\eta=X_1X_2$. Then both $\xi$ and $\eta$ belong to the second Wiener chaos of $W$ (write 
	\[ \frac{1}{2} (X_1^2-1) + \frac{1}{2}(X_2^2-1) - \left[\left(\frac{X_1-X_2}{\sqrt{2}}\right)^2-1\right] = X_1X_2\]
for $\eta$), $\mathbb{E}\xi=\mathbb{E}\eta=0$, and $\mathbb{E}\xi^2=\mathbb{E}\eta^2=1$. Then we have that 
	\[ |\mathbb{E}(\xi^2-\mathbb{E}\xi^2)(\eta^2-\mathbb{E}\eta^2)| = 4\]
but $\mathbb{E}\xi\eta = 0$, hence $\rho_{\xi,\eta}=0$, and an inequality of the type \eqref{eq:Gebelein} does not hold. 
\end{counterexample}

Altogether, it turns out that such a strong tool as Gebelein's inequality does not hold in necessary power already in the second Wiener chaos and therefore, in order to prove \autoref{main_thm}, we use orthogonality of Wiener chaoses and tensor calculus instead. In particular, we work directly with the Besov-Orlicz norm; we initially use the generalized product formula for multiple integrals (given in \autoref{exp_thm}) to obtain the Wiener chaos expansion of
	\[ Y_{\ell,\delta}^\ell=\|G(\cdot+\delta)-G(\cdot)\|_{L^\ell(0,T-\delta)}^\ell\]
for an even integer $\ell$ in terms of tensor cancellations of the kernel $A_t$, i.e.
	\[ Y_{\ell,\delta}^\ell = \underbrace{\mathbb{E}Y_{\ell,\delta}^\ell}_{\in\mathscr{H}_0} + \underbrace{\xi_2(\delta)}_{\in\mathscr{H}_2} +\underbrace{\xi_4(\delta)}_{\in\mathscr{H}_4} + \ldots + \underbrace{\xi_{\ell n}(\delta)}_{\in\mathscr{H}_{\ell n}}\]
where $\mathscr{H}_k$ denotes the $k$\textsuperscript{th} Wiener chaos of the isonormal process $W$. Subsequently, by using several results regarding the cancellation operator (and, in particular, the key \autoref{mult_est_2} that is used instead of Gebelein's inequality), we show\footnote{More precisely, we show that the variances $\mathbb{E}[\xi_{k}(\delta)]^2$, $k=2,4,\ldots, \ell n$, are negligible so that the mean-square distance $\mathbb{E}(Y_{\ell,\delta}^\ell - \mathbb{E}Y_{\ell,\delta}^\ell)^2$ is small by orthogonality of Wiener chaoses.} that 
	\[ |\xi_{2}(\delta)|+|\xi_4(\delta)| + \ldots + |\xi_{\ell n}(\delta)| \ll \mathbb{E}Y_{\ell,\delta}^\ell \quad\mbox{as}\quad\delta\rightarrow 0+\]
so that \[Y_{\ell,\delta}^\ell \sim \mathbb{E}Y_{\ell,\delta}^\ell \quad \mbox{as}\quad \delta\to 0+.\] From this, we obtain that 
\[Y_{\ell,\delta}\lesssim \delta^\alpha\ell^\frac{n}{2}\] and the result follows.

\subsection{Organisation of the article}

In \autoref{sec:prelim}, the definitions of Besov, Orlicz, and Besov-Orlicz spaces are recalled. In \autoref{sec:main}, the main results of the article are collected. In particular, we state the main \autoref{main_thm} where, along some moment estimates, we give sufficient conditions for the considered process to have paths in the Besov-Orlicz space $B_{\varPhi_{2/n},\infty}^\alpha(0,T)$ and sufficient conditions for the process not to have paths in any of the Besov spaces $B_{p,q}^\alpha(0,T)$ for $p\in [1,\infty]$ and $q\in [1,\infty)$. In addition, we also give \autoref{cor:mod_holder} of the theorem to assess modular H\"older continuity of the process. Subsequently, we discuss the assumptions of the theorem and compare it to the known results in \autoref{rem:G3'}. In there, we also note that some alternative assumptions, that correspond to those in \cite[Theorem 5.1]{Ver09}, can be considered. The section is concluded by \autoref{thm:main_thm_2} which refines \autoref{main_thm} and in which pathwise asymptotics of the integral increments of the considered process is treated. In \autoref{sec:examples}, the conditions are verified for fractionally filtered Hermite processes with a product kernel. In \autoref{sec:tensor_calc}, we review elements of tensor calculus, provide several motivating examples and the necessary technical tools for the proofs of \autoref{main_thm} and \autoref{thm:main_thm_2} (\autoref{lem_2}, \autoref{mult_est}, and \autoref{mult_est_2} regarding the properties of the cancellation operator and \autoref{exp_thm} that contains the generalized product formula for multiple integrals). The proofs of \autoref{main_thm}, \autoref{main_thm} with alternative assumptions, and \autoref{thm:main_thm_2} are given at the end of the article in \autoref{sec:proof}, \autoref{sec:app_A}, and \autoref{sec:proof_2}, respectively.

\section{Preliminaries: Besov, Orlicz, and Besov-Orlicz spaces}
\label{sec:prelim}

We begin by recalling some facts about Besov, Orlicz, and Besov-Orlicz spaces. For a thorough exposition on Besov spaces, we refer the reader to, e.g., \cite{Trie83}. Orlicz spaces are covered in, e.g., \cite{RaoRen91,Zaa83} and Besov-Orlitz spaces in, e.g., \cite{OndVer20,PickSick93}. Let $I\subseteq [0,\infty)$ be a bounded interval and for $h\in\mathbb{R}$, we denote $I(h)=\{s\in I: s+h\in I\}$. For $s\in (0,1)$ and $p,q\in [1,\infty)$ (with the usual modifications for $p=\infty$ or $q=\infty$), the Besov space $B_{p,q}^s(I)$ is defined as the linear space
	\[ B_{p,q}^s(I)=\{f\in L^p(I): (f)_{B_{p,q}^s(I)}<\infty\}\]
where 
	\[ (f)_{B_{p,q}^s(I)} = \left(\int_{0}^{{\infty}} [t^{-s}\sup_{|h|\leq t}\|f(\cdot + h)-f(\cdot)\|_{L^p(I(h))}]^q\frac{\d{t}}{t}\right)^\frac{1}{q}\]
and the space $B_{p,q}^s(I)$ is a Banach space when endowed with the norm
	\[ \|f\|_{B_{p,q}^s(I)}= \|f\|_{L^p(I)} + (f)_{B_{p,q}^s(I)}.\]
Besov spaces can be generalized as follows. A function $\mathcal{N}:[0,\infty)\rightarrow\mathbb{R}$ is called a \textit{Young function} if it is non-negative, non-decreasing, continuous, convex, and satisfies $\mathcal{N}(0)=0$ and $\mathcal{N}(\infty-)=\infty$. For a Young function $\mathcal{N}$ and a measure space $(D,\mathcal{D},\mu)$, where $\mu$ is a $\sigma$-finite measure, the \textit{Orlicz space (with the Young function $\mathcal{N}$)} $L^\mathcal{N}(D)$ is defined as the linear space 
	\[ L^\mathcal{N}(D)= \{f\in L^0(D): \|f\|_{L^{\mathcal{N}}(D)}<\infty\}\]
where 
	\[ \|f\|_{L^\mathcal{N}(D)} = \inf\left\{\lambda\geq 0: \int_{D} \mathcal{N}\left(\frac{|f(x)|}{\lambda}\right)\mu(\d{x}) \leq 1\right\}\]
is the so-called \textit{Luxemburg norm}. Endowed with this norm, $L^\mathcal{N}(D)$ is a Banach space. Finally, for $s\in (0,1)$, a Young function $\mathcal{N}$, and $q\in[1,\infty)$ (with the usual modification for $q=\infty$), the \textit{Besov-Orlicz} space $B_{\mathcal{N},q}^s(I)$ is defined as the linear space 
	\[ B_{\mathcal{N},q}^s(I)= \{f\in L^\mathcal{N}(I): (f)_{B_{\mathcal{N},q}^s(I)}<\infty\}\]
where 
	\[ (f)_{B_{\mathcal{N},q}^s(I)} = \left(\int_0^\infty [t^{-s}\sup_{|h|\leq t}\|f(\cdot + h)-f(\cdot)\|_{L^\mathcal{N}(I(h))}]^q\frac{\d{t}}{t}\right)^\frac{1}{q}.\]
The space $B_{\mathcal{N},q}^s(I)$ is a Banach space when endowed with the norm
	\[ \|f\|_{B_{\mathcal{N},q}^s(I)} = \|f\|_{L^{\mathcal{N}}(I)} + (f)_{B_{\mathcal{N},q}^s(I)}.\]
\begin{remark} There are equivalent (semi)norms that may be more convenient in certain problems. In particular, for $s\in (0,1)$ and $p,q\in [1,\infty)$ (with the usual modifications for $p=\infty$ or $q=\infty$), the seminorm $(f)_{B_{p,q}^s(I)}$ is equivalent to 
	\[ [f]_{B_{p,q}^s(I)} = \left(\sum_{j\geq 0} 2^{jsq}\|f(\cdot + 2^{-j})-f(\cdot)\|_{L^p(I(2^{-j}))}^q\right)^\frac{1}{q}\]
by dyadic approximation; see, e.g., \cite[Corollary 3.b.9]{Kon86}. Moreover, in the present article, particular attention will be given to the \textit{exponential} Orlicz and Besov-Orlicz spaces $L^{\varPhi_\beta}(D)$ and $B_{\varPhi_\beta,\infty}^s(I)$ for $s\in (0,1)$ where $\varPhi_\beta$, $\beta>0$, is a Young function that satisfies $\varPhi_\beta(x)=\mathrm{e}^{x^\beta}-1$ for $x\in [\tau_\beta,\infty)$ with some $\tau_\beta\geq 0$. In this case, the (semi)norms $\|f\|_{L^{\varPhi_\beta}(D)}$ and $(f)_{B_{\varPhi_\beta,\infty}^s(I)}$ are equivalent to 
	\[ \vertiii{f}_{L^{\varPhi_{\beta}}(D)} = \sup_{p\geq 1} p^{-1/\beta}\|f\|_{L^p(D)}\]
and
	\[ [f]_{B_{\varPhi_\beta,\infty}^s(I)} =\sup_{j\geq 1}2^{js}\vertiii{f(\cdot + 2^{-j})-f(\cdot)}_{L^{\varPhi_\beta}(I(2^{-j}))},\]
respectively, by, e.g., \cite[Proposition 2.3]{OndVer20}. It follows that the norm $\|f\|_{B_{\varPhi_\beta,\infty}^s(I)}$ is equivalent to 
	\[ \vertiii{f}_{B_{\varPhi_\beta,\infty}^s(I)} = \vertiii{f}_{L^{\varPhi_\beta}(I)} + [f]_{B_{\varPhi_{\beta},\infty}^s(I)}.\]
\end{remark}
	
\section{Main results: Path regularity}
\label{sec:main}

The main results of the article are collected in this section. Let $H$ be a real separable Hilbert space, $W=\{W(h)\}_{h\in H}$ be an $H$-isonormal Gaussian process, i.e. a centered Gaussian process with the covariance 
	\[
		\mathbb E\,[W(h_1)W(h_2)]=\langle h_1,h_2\rangle_{H},\qquad h_1,h_2\in H,
	\]
defined on a complete probability space $(\Omega, \mathcal{F},\mathbb{P})$, and assume that the $\sigma$-algebra $\mathcal{F}$ is generated by process $W$. For $k\in\mathbb{N}_0$ we denote by
	\[
		W_k:H^{\otimes_2^k}\to L^2(\Omega)
	\]
the multiple divergence operator of order $k$ defined via the duality
	\[
		\mathbb E\,\langle D^kX,A\rangle_{H^{\otimes_2^k}}=\mathbb E\,[XW_k(A)],\qquad A\in H^{\otimes_2^k},
	\]
for all $X\in\mathbb D^{k,2}$ where $D^k$ denotes the $k$\textsuperscript{th} Malliavin derivative and $\mathbb D^{k,2}$ its domain; see, e.g., \cite{NouPec12} or \cite{Nua06} for details. Let us fix $T>0$, $n\in\mathbb{N}$, and $\alpha\in (0,1)$ for the rest of this section. Let $G=\{G(t)\}_{t\in [0,T]}$ be a jointly measurable process in the $n$\textsuperscript{th} Wiener chaos of $W$ represented by
\begin{equation}
\label{eq:G}
	G(t)=W_n(A_t),\qquad t\in[0,T],
\end{equation}
where $A_t\in H^{\otimes_2^n}$ is (the unique) symmetric tensor and $t\mapsto A_t$ is a Bochner measurable function from $[0,T]$ to $H^{\otimes_2^n}$. We define
\[
A_{x,s}=A_{x+s}-A_x
\]
for $x,s\geq 0$, $x+s\leq T$, and 
\[ C_{s,t}(x,y) = s^{-2\alpha}t^{-2\alpha} \|A_{x,s}\otimes _1A_{y,t}\|_{H^{\otimes_2^{2n-2}}}^2 + s^{-\alpha}t^{-\alpha} \sum_{j=2}^n \|A_{x,s}\otimes_jA_{y,t}\|_{H^{\otimes_2^{2(n-j)}}}\]
for $x,y\geq 0$, $s,t>0$ that satisfy $x+s\leq T$, $y+t\leq T$. Here, $\otimes_j$, $j=0,\ldots, n$, is the cancellation operator defined by
\begin{align*}
(h_1\otimes\dots\otimes h_n)\otimes_0(k_1\otimes\dots\otimes k_n)&=h_1\otimes\dots\otimes h_n\otimes k_1\otimes\dots\otimes k_n,
\\
(h_1\otimes\dots\otimes h_n)\otimes_j(k_1\otimes\dots\otimes k_n)&=\langle h_1,k_1\rangle_{H}\dots\langle h_j,k_j\rangle_{H}\,h_{j+1}\otimes\dots\otimes h_n\otimes k_{j+1}\otimes\dots\otimes k_n,
\\
(h_1\otimes\dots\otimes h_n)\otimes_n(k_1\otimes\dots\otimes k_n)&=\langle h_1,k_1\rangle_{H}\dots\langle h_n,k_n\rangle_{H}
\end{align*}
for $h_i,k_i\in H$, $i=1,\ldots, n$. We also define
\[ F (s,t) = \int_0^{T-t}\int_0^{T-s}C_{s,t}(x,y)\d{x}\d{y}\]
for $s,t\in (0,T)$.

\medskip

In the rest of the article, various combinations of assumptions \ref{ass:G1} - \ref{ass:G4} below are considered and these assumptions are formulated now. 

\begin{assumption} 
There exist $\kappa\in [1,\infty)$, $\kappa^\prime\in (0,\infty)$, and $\varepsilon\in (0,1)$ such that
\begin{enumerate}[label=(G\arabic*)]
	\item\label{ass:G1} $\|A_{x,s}\|_{H^{\otimes_2^n}}\le\kappa s^\alpha$ for every $x\ge 0$, $s>0$, $x+s\le T$,
	\item\label{ass:G2} $ \liminf_{s\to 0}\inf_{x\in [0,T)} s^{-\alpha}\|A_{x,s}\|_{H^{\otimes_2^n}}\geq \kappa^\prime$, 
	\item\label{ass:G3} $ \sum_{j=j_0}^\infty F(2^{-j},2^{-j})<\infty$ where $j_0=\min\{j: 2^{-j}< T\}$,
	\item\label{ass:G4} $F(s,t)\leq \kappa s^{\varepsilon}t^{\varepsilon}$ for every $s,t\in (0,T)$.
\end{enumerate}
\end{assumption}

\begin{remark}
Note that if the process $G$ satisfies condition \ref{ass:G1}, then it has paths in the Besov space 
	\[ B_{p,q}^r(0,T)\]
for every $p,q\in [1,\infty]$ and $r\in (0,\alpha)$. Indeed, by \cite[Theorem 2.7.2]{NouPec12} and Kolmogorov's continuity criterion, $G$ has paths in the H\"older space $C^r([0,T])$ for every $r\in (0,\alpha)$ and this is equivalent to the claim by the embedding of Besov spaces from \cite[Theorem 3.3.1 and Proposition 3.2.4]{Trie83}.
\end{remark}

Consider now Young functions $\varPhi_\beta$ for $\beta>0$ that satisfy $\varPhi_\beta(x)=e^{x^\beta}-1$ for $x\in [\tau_\beta,\infty)$ with some $\tau_\beta\ge 0$. The main result of the paper follows. Its proof is postponed to \autoref{sec:proof}.

\begin{theorem}
\label{main_thm} 
If process $G$ satisfies \ref{ass:G1} and \ref{ass:G3}, then 
	\[ \|G\|_{B_{\varPhi_{2/n},\infty}^\alpha(0,T)} \in L^{\varPhi_{2/n}}(\Omega)\]
and, in particular, it has paths in the Besov-Orlicz space $B^\alpha_{\varPhi_{2/n},\infty}(0,T)$
almost surely. If, additionally, process $G$ satisfies \ref{ass:G2}, then its paths do not belong to the Besov space $B^{\alpha}_{p,q}(0,T)$ for any $p\in[1,\infty]$ and $q\in [1,\infty)$ almost surely.
\end{theorem}

By the results contained in \cite[Section 2.4.2]{OndVer20}, the following corollary is obtained.

\begin{corollary}
\label{cor:mod_holder} If $G$ satisfies \ref{ass:G1} and \ref{ass:G3}, then it has paths in the modular H\"older space
\[
C^{r^\alpha|\log r|^\frac{n}{2}}([0,T])
\]
almost surely; i.e. there exists an almost surely finite and positive random variable $\bm{c}$ such that
\[
|G(s)-G(t)|\le \bm{c}\,|s-t|^\alpha|\log |s-t||^\frac{n}{2},\qquad s,t\in[0,T], \,0<|s-t|<1/2,\quad\text{a.s.}
\]
\end{corollary}

\begin{remark}
\label{rem:HSSI}
Note that every $\alpha$-self-similar process $G$ of the form \eqref{eq:G} that has stationary increments satisfies conditions \ref{ass:G1} and \ref{ass:G2}. Indeed, in this case, we have 
	\[
		\|A_{x,s}\|_{H^{\otimes_2^n}}^2  = \frac{1}{n!}\E[G(x+s)-G(x)]^2 = \frac{1}{n!} \E[G(s)-G(0)]^2= \frac{1}{n!}\E[G(1)]^2 s^{2\alpha}
	\]
for every $x\geq 0$, $s > 0$, $x+s\leq T$ where stationarity of increments of $G$, the fact that $G(0)=0$ holds almost surely, and self-similarity of $G$ were used successively. 
\end{remark}

\begin{remark}
\label{rem:G3'} Note also that there is a discrete version of condition \ref{ass:G3} in \autoref{main_thm}. In particular, \autoref{main_thm} remains valid if condition \ref{ass:G3} is replaced by
\begin{enumerate}[label=(G3$^\prime$)]
\item\label{ass:G3'} There exists a function $K:[0,T]^2\to\mathbb [0,\infty)$ such that \[C_{s,s}(x,y)\leq K(s, |x-y|)\] holds for $x,y\geq 0$, $s>0$ that satisfy $x+s\leq T$, $y+s\leq T$, $(x,x+s)\cap (y,y+s)=\emptyset$, and 
		\[ \sum_{j=1}^\infty \delta_j^2 \sum_{\substack{m,m^\prime=1\\ m\neq m^\prime}}^{J_{\delta_j}} K(\delta_j,|m-m^\prime|\delta_j)<\infty\]
	holds with 
		\[\delta_j=T2^{-j}\quad\mbox{and}\quad J_{\delta_j}=2^j-1.\]
\end{enumerate}
The modified proof can be found in \autoref{sec:app_A}. In the Gaussian case, the above condition \ref{ass:G3'} is implied by the conditions in \cite[Theorem 5.1]{Ver09} and, consequently, \autoref{main_thm} provides an extension of that result in the scalar case to the general setting.
\end{remark}

Clearly, to verify condition \ref{ass:G3} it suffices to find $\varepsilon\in (0,1)$ such that $F(s,s)\leq \kappa s^{2\varepsilon}$ holds for every $s\in (0,T)$. If this condition is slightly strengthened, we obtain the exact pathwise asymptotics of the integral increments of $G$. The proof is postponed to \autoref{sec:proof_2}.

\begin{notation}
We write
	\[f(x)\in\Theta(g(x))\quad\mbox{as}\quad x\to a+\qquad\mbox{if}\qquad 0<\liminf_{x\rightarrow a+} \frac{f(x)}{g(x)} \leq \limsup_{x\rightarrow a+}\frac{f(x)}{g(x)} <\infty.\]
\end{notation}

\begin{theorem}
\label{thm:main_thm_2}
If process $G$ satisfies \ref{ass:G1}, \ref{ass:G2}, and \ref{ass:G4}, then, in addition to the assertions of \autoref{main_thm} being true, we have that
	\[ \|G(\cdot+s)-G(\cdot)\|_{L^p(0,T-s)} \in \Theta(s^\alpha)\quad \mbox{as}\quad s\to 0+\quad\mbox{a.s.}\]
holds for every $p\in [1,\infty)$ and 
	\[ \|G(\cdot+s)-G(\cdot)\|_{L^{\varPhi_{2/n}}(0,T-s)} \in \Theta(s^\alpha)\quad \mbox{as}\quad s\to 0+\quad\mbox{a.s.}\]
\end{theorem}

\begin{remark}
Although condition \ref{ass:G4} is by no means sharp it suffices for the demonstration of the method in \autoref{thm:main_thm_2} and it covers all our examples. Therefore, only this criterion is given here for simplicity of the exposition.
\end{remark}

\begin{notation}
Throughout the article, we write $A\lesssim B$ (and $A\gtrsim B$) whenever there exists a finite positive constant $C$ such that $A\leq CB$ (and $A\geq CB$) whose precise value is not important. This constant can also change from line to line. We also write $A\lesssim_K B$ (and $A\gtrsim_K B$) to indicate the dependence of constant $C$ on a different quantity $K$. If $C$ is random, we indicate this by writing $A\leom B$ (and $A\geom B$).
\end{notation}

\section{Example: Fractionally filtered generalized Hermite processes}
\label{sec:examples}

\subsection{Application to fractionally filtered generalized Hermite processes} In this section, the results are applied to a specific class of stochastic processes that contains some well-studied examples such as fractional Brownian motions or Hermite processes. Let $H=L^2(\mathbb R)$ and let
	\begin{equation}
	\label{eq:A_Hermite}		 
		 A_t(x_1,\ldots,x_n) = \int_\mathbb{R}k_t^{\beta_1}(u)\prod_{i=1}^n\phi^{\beta_2}(u-x_i)\d{u}, \qquad x_1,\ldots,x_n\in \mathbb R,\quad  t\in [0,T],
	\end{equation}
for the function $k_t^{\beta_1}:\mathbb{R}\rightarrow\mathbb{R}$ defined by 
	\begin{equation}
	\label{eq:k_t} k_t^{\beta_1}(u)=\begin{cases}
						\bm{1}_{(0,t]}(u), & \quad \beta_1=0,\\
						\frac{1}{\beta_1}[(t-u)_+^{\beta_1}-(-u)_+^{\beta_1}], & \quad \beta_1\neq 0,
						\end{cases}
	\end{equation}
where $(u)_+=\max\{0,u\}$, and for some measurable function $\phi^{\beta_2}: \mathbb R\rightarrow\mathbb R$ that satisfies
	\begin{enumerate}[label=(H1)]
	\item\label{ass:H1}
	\begin{itemize}
	\item $\int_{\mathbb R} |\phi^{\beta_2}(x)\phi^{\beta_2}(x+u)|\d{x}\leq \kappa u^{\beta_2-1}$, $u>0,$ with $\kappa\in (0,\infty)$ and $\beta_2\in (1-\frac{1}{n}, 1)$.
	\end{itemize}
	\end{enumerate}
	
We also assume that

	\begin{enumerate}[label=(H2)]
	\item\label{ass:H11}
	\begin{itemize}
	\item $0<\alpha=\beta_1 + \frac{n}{2}(\beta_2-1) + 1<1$.
	\end{itemize}
	\end{enumerate}
	
In this situation, the following result is obtained.

\begin{corollary}
\label{cor:Hermite_1}
Process $G$ defined by \eqref{eq:G} with kernel $A$ defined by \eqref{eq:A_Hermite} for which assumptions \ref{ass:H1} and \ref{ass:H11} hold satisfies conditions \ref{ass:G1} and \ref{ass:G4}, and, consequently, it holds that 
	\[ \|G\|_{B_{\varPhi_{2/n},\infty}^\alpha(0,T)} \in L^{\varPhi_{2/n}}(\Omega).\]
\end{corollary}

\begin{proof}
We only treat the case of $\beta_1\neq 0$; the case $\beta_1=0$ follows by similar arguments. Denote
	\[K(u)=\int_{\mathbb{R}} \phi^{\beta_2}(x)\phi^{\beta_2}(x+u)\d{x}.\]
Then $K$ is a symmetric locally bounded function on $\mathbb{R}\setminus\{0\}$ and with the notation introduced in \autoref{sec:main}, we have that 
	\[
		 \|A_{x,s}\|_{L^2(\mathbb R^n)}^2  =\int_{\mathbb{R}^2} k_s^{\beta_1}(u)k_{s}^{\beta_1}(v)K(|u-v|)^n\d{u}\d{v}
	\]
holds for every $x \geq 0$ and $s>0$ such that $x+s\leq T$. As there is the inequality
	\begin{equation}
	\label{eq:est_G1}
		 \int_{\mathbb{R}^2} |k_s^{\beta_1}(u)||k_{s}^{\beta_1}(v)| |K(|u-v|)|^n\d{u}\d{v} \lesssim s^{2\alpha}
	\end{equation}
by assumptions \ref{ass:H1} and \ref{ass:H11}, it follows that condition \ref{ass:G1} is satisfied. In order to verify condition \ref{ass:G4} notice first that 
	\[ |k_s^{\beta_1}(-v)| \approx \begin{cases}
														(s+v)^{\beta_1} , & \quad v\in (-s,0],\\
														\begin{cases}
															s^{\beta_1}, & \quad \beta_1>0,\\
															v^{\beta_1}, & \quad \beta_1<0,
														\end{cases} & \quad v\in (0,s],\\
														sv^{\beta_1-1}, & \quad v\in [s,\infty),
													\end{cases}\]
where the approximation constants do not depend on either $s$ or $v$. Now, denote
	\[|x|^{\gamma_1,\gamma_2} = \bm{1}_{[|x|<1]}|x|^{\gamma_1} + \bm{1}_{[|x|\geq 1]}|x|^{\gamma_2}, \qquad x\in\mathbb{R}\setminus\{0\},\]
for $\gamma_1,\gamma_2\in\mathbb R$. It follows from the above approximation that
	\[ \int_{\mathbb{R}} |k_s^{\beta_1}(-v)||u-v|^{0,n(\beta_2-1)}\d{v} \approx s^{1+\min\{\beta_1,0\}} |u|^{0,n(\beta_2-1)+\max\{\beta_1,0\}}\]
from which we obtain the estimate
	\begin{equation}
	\label{eq:est_Q}
		Q(s,t) = \int_{\mathbb{R}^2} |k_t^{\beta_1}(u)||k_s^{\beta_1}(v)| |u-v|^{0,n(\beta_2-1)}\d{v}\d{u} \lesssim s^{1+\min\{\beta_1,0\}}t^{1+\min\{\beta_1,0\}}, \qquad s,t\in [0,T].
	\end{equation}
We are now in position to verify condition \ref{ass:G4}. Observe that there is the estimate 
	\[s^{-2\alpha}t^{-2\alpha}\|A_{x,s}\otimes_1A_{y,t}\|_{L^2(\mathbb{R}^{2n-2})}^2 \lesssim s^{-\alpha} t^{-\alpha} \|A_{x,s}\otimes_1A_{y,t}\|_{L^2(\mathbb{R}^{2n-2})}\]
for $x,y\geq 0$, $s,t>0$ such that $x+s\leq T$, $y+t\leq T$ by \ref{ass:G1} and the forthcoming \autoref{mult_est}. This means that it suffices to find $\varepsilon\in (0,1)$ such that 
	\[\int_0^{T-t}\int_0^{T-s} \sum_{j=1}^n \|A_{x,s}\otimes_j A_{y,t}\|_{L^2(\mathbb{R}^{2(n-j)})}\d{x}\d{y} \lesssim s^{\alpha+\varepsilon}t^{\alpha+\varepsilon}\]
for $s,t\in (0,T)$. We have that
	\begin{multline*}
		 \|A_{x,s}\otimes_jA_{y,t}\|_{L^2(\mathbb{R}^{2n-2j})}^2 = \int_{\mathbb{R}^4} k_s^{\beta_1}(r_1)k_t^{\beta_1}(r_2)k_s^{\beta_1}(r_3)k_t^{\beta_1}(r_4)\\
		 \cdot K(|r_1-r_3|)^{n-j}K(|r_2-r_4|)^{n-j} K(|r_1-r_2+x-y|)^j K(|r_3-r_4+x-y|)^j\d{r}
	\end{multline*}
and, consequently, that
	\begin{align*}
	 \|A_{x,s}\otimes_jA_{y,t}\|_{L^2(\mathbb{R}^{2n-2j})}^2 & \\
	 & \hspace{-3.2cm} \leq \int_{\mathbb{R}^4} \left( |k_s^{\beta_1}(r_1)k_t^{\beta_1}(r_2)k_s^{\beta_1}(r_3)k_t^{\beta_1}(r_4)| |K(|r_1-r_3|)|^n|K(|r_2-r_4|)|^n\right)^{1-\frac{j}{n}}\\
	 & \hspace{-2.7cm} \cdot \left(| k_s^{\beta_1}(r_1)k_t^{\beta_1}(r_2)k_s^{\beta_1}(r_3)k_t^{\beta_1}(r_4)| |K(|r_1-r_2+x-y|)|^n |K(|r_3-r_4+x-y|)|^n\right)^\frac{j}{n}\d{r}\\
	 & \hspace{-3.2cm} \leq \left(\int_{\mathbb{R}^4} |k_s^{\beta_1}(r_1)k_t^{\beta_1}(r_2)k_s^{\beta_1}(r_3)k_t^{\beta_1}(r_4)||K(|r_1-r_3|)|^n|K(|r_2-r_4|)|^n\d{r}\right)^{1-\frac{j}{n}}\\
	 & \hspace{-2.7cm} \cdot\left(\int_{\mathbb{R}^4}|k_s^{\beta_1}(r_1)k_t^{\beta_1}(r_2)k_s^{\beta_1}(r_3)k_t^{\beta_1}(r_4)||K(|r_1-r_2+x-y|)|^n|K(|r_3-r_4+x-y|)|^n\d{r}\right)^\frac{j}{n}\\
	 & \hspace{-3.2cm} = \left(\int_{\mathbb{R}^2} |k_{s}^{\beta_1}(u)||k_{s}^{\beta_1}(v)| |K(|u-v|)|^n\d{u}\d{v}\right)^{1-\frac{j}{n}}  \left(\int_{\mathbb{R}^2} |k_{t}^{\beta_1}(u)||k_{t}^{\beta_1}(v)| |K(|u-v|)|^n\d{u}\d{v}\right)^{1-\frac{j}{n}}\\
	 & \hspace{-2.7cm} \cdot\left(\int_{\mathbb{R}^2}|k_{s}^{\beta_1}(u)||k_{t}^{\beta_1}(v)||K(|u-v+x-y|)|^n\d{u}\d{v}\right)^\frac{2j}{n}
	\end{align*}
hold for every $x,y\geq 0$ and $s,t>0$ such that $x+s\leq T$, $y+t\leq T$, and $j\in\{1,\ldots,n\}$ by using H\"older's inequality to obtain the second estimate. Hence, the successive use of the above estimate, inequality \eqref{eq:est_G1}, Jensen's inequality, assumption \ref{ass:H1}, and estimate \eqref{eq:est_Q} yields 
	\begin{align*}
		\int_0^{T-t}\int_0^{T-s} \|A_{x,s}\otimes_jA_{y,t}\|_{L^2(\mathbb{R}^{2n-2j})}\d{x}\d{y} & \\
		& \hspace{-5cm} \lesssim s^{\alpha\left(1-\frac{j}{n}\right)}t^{\alpha\left(1-\frac{j}{n}\right)}\left(\int_0^T\int_0^T \int_{\mathbb{R}^2}|k_s^{\beta_1}(u)||k_t^{\beta_1}(v)| |K(|u-v+x-y|)|^n\d{u}\d{v}\d{x}\d{y}\right)^{\frac{j}{n}} \\
		& \hspace{-5cm} \lesssim s^{\alpha\left(1-\frac{j}{n}\right)}t^{\alpha\left(1-\frac{j}{n}\right)}[Q(s,t)]^\frac{j}{n}\\
		& \hspace{-5cm} \lesssim s^{\alpha + [1+\min\{\beta_1,0\} -\alpha]\frac{j}{n}}t^{\alpha + [1+\min\{\beta_1,0\} -\alpha]\frac{j}{n}}
	\end{align*}
for $s,t\in (0,T)$ and $j\in \{1,\ldots, n\}$. As $1+\min\{\beta_1,0\}-\alpha>0$ by \ref{ass:H1} and \ref{ass:H11}, condition \ref{ass:G4} is verified and the claim of the corollary is obtained by appealing to \autoref{main_thm}.
\end{proof}

\subsection{Application to fractional Brownian motions and Hermite processes}

The fact that sample paths of the fractional Brownian motion with Hurst parameter $\alpha\in (0,1)$ belong to the Besov-Orlicz space $B_{\varPhi_2,\infty}^\alpha(0,1)$ and do not belong to Besov space $B_{p,q}^\alpha(0,1)$ for any $p\in [1,\infty]$ and $q\in [1,\infty)$ has already been established in \cite[Corollary 5.3]{Ver09}. In addition to this, we improve the result by obtaining estimates on the moments of its norm and the asymptotics of its integral increments. Let $\beta_1$ and $\beta_2$ be such that 
	\[ 0<\beta_1+\frac{1}{2}(\beta_2-1)+1<1\quad\mbox{and}\quad 0<\beta_2<1\]
for the purposes of the next corollary. Then the process $z_1^{\beta_1,\beta_2}$ defined by
	\[ z_1^{\beta_1,\beta_2}(t) = c_{\beta_1,\beta_2}^{(1)}\int_{\mathbb{R}}\left[\int_{\mathbb{R}} k_t^{\beta_1}(u)(u-x)_+^{\frac{\beta_2}{2}-1}\d{u}\right]\d{W}_x, \qquad t\in [0,1],\]
where $\int_{\mathbb{R}}(\cdots)\d{W}_x$ is a Wiener-It\^o integral with respect to the Wiener process $\{W_t\}_{t\in\mathbb{R}}$ and $c_{\beta_1,\beta_2}^{(1)}$ is a suitable normalising constant (cf., e.g., \cite{BaiTaq14} or \cite{MaeTud07}) is a fractional Brownian motion with Hurst parameter $\alpha=\beta_1+\frac{1}{2}(\beta_2-1) +1$ and we have
\begin{corollary}
It holds for process $z_{1}^{\beta_1,\beta_2}$ that
	\[ \|{z_1^{\beta_1,\beta_2}}\|_{B_{\varPhi_{2},\infty}^\alpha(0,T)} \in L^{\varPhi_{2}}(\Omega).\]
Moreover, it holds that
	\[ \|z_1^{\beta_1,\beta_2}(\cdot+s)-z_1^{\beta_1,\beta_2}(s)\|_{L^p(0,T-s)}\in \Theta(s^\alpha)\quad\mbox{as}\quad s\to 0+\]
for every $p\in [1,\infty)$ almost surely and
	\[ \|z_1^{\beta_1,\beta_2}(\cdot+s)-z_1^{\beta_1,\beta_2}(s)\|_{L^{\varPhi_{2}}(0,T-s)}\in \Theta(s^\alpha)\quad\mbox{as}\quad s\to 0+\]
almost surely.
\end{corollary}
\begin{proof}
Conditions \ref{ass:G1} and \ref{ass:G2} are verified by appealing to \autoref{rem:HSSI} (process $z_1^{\beta_1,\beta_2}$ is $\alpha$-self-similar and has stationary increments) and condition \ref{ass:G4} holds by \autoref{cor:Hermite_1}. The claim then follows by \autoref{thm:main_thm_2}.
\end{proof}
The results are applied to Hermite processes now. Let $n\in\mathbb{N}$, $n\geq 2$, and let $\alpha\in (1/2,1)$ for the purposes of the next corollary. The process $\{z_n^{\alpha}(t)\}_{t\in [0,1]}$ defined by
	\[ z_n^{\alpha}(t) = c_{\alpha,n}^{(2)}\int_{\mathbb{R}^n}\left[\int_0^t \prod_{i=1}^n(u-x_i)_+^{-\frac{1}{2}-\frac{1-\alpha}{n}}\d{u}\right]\d{W}_{\bm{x}}, \qquad t\in [0,1],\]
where $\int_{\mathbb{R}^n}(\cdots)\d{W}_{\bm{x}}$ is the multiple Wiener-It\^o integral with respect to the Wiener process $\{W_t\}_{t\in \mathbb{R}}$ of order $n$ and $c_{\alpha,n}^{(2)}$ is a suitable normalizing constant (cf., e.g., \cite{BaiTaq14} or \cite{MaeTud07}), is a Hermite process of order $n$ and Hurst parameter $\alpha$. We have
\begin{corollary}
 It holds for the Hermite process $z_n^\alpha$ that
	\[ \|z_n^\alpha\|_{B_{\varPhi_{2/n},\infty}^\alpha(0,T)} \in L^{\varPhi_{2/n}}(\Omega)\] and it does not have sample paths in the Besov space $B_{p,q}^\alpha(0,1)$ with $p\in [1,\infty]$ and $q\in [1,\infty)$ almost surely. Moreover, we have that 
	\[ \|z_n^{\alpha}(\cdot + s)-z_n^\alpha(s)\|_{L^p(0,T-s)}\in \Theta(s^\alpha) \quad\mbox{as}\quad s\to 0+\]
holds for every $p\in [1,\infty)$ almost surely and 
	\[ \|z_n^{\alpha}(\cdot + s)-z_n^\alpha(s)\|_{L^{\varPhi_{2/n}}(0,T-s)}\in \Theta(s^\alpha)\quad\mbox{as}\quad s\to 0+\]
holds almost surely.
\end{corollary}

\begin{proof}
Conditions \ref{ass:G1} and \ref{ass:G2} are verified by appealing to \autoref{rem:HSSI} (process $z_n^{\alpha}$ is $\alpha$-self-similar and has stationary increments) and condition \ref{ass:G4} holds by \autoref{cor:Hermite_1}. The claim then follows by \autoref{thm:main_thm_2}.
\end{proof}
\begin{remark}
In \cite{BaiTaq14}, the following family of processes is considered. Define the kernel $A$ by
	\begin{equation}
	\label{eq:A_Hermite_general}
	 A_t(x_1,\ldots,x_n)=\int_{x_1\vee\ldots\vee x_n}^\infty k_t^{\gamma_1}(u)\phi^{\gamma_2}(u-x_1,\ldots, u-x_n)\d{u}, \qquad x_1,\ldots, x_n\in\mathbb{R},\quad t\in [0,T],
	 \end{equation}
where $k_t^{\gamma_1}$ is the function defined by formula \eqref{eq:k_t} and where $\phi^{\gamma_2}:\mathbb{R}^n_+\rightarrow \mathbb{R}$ is a non-zero function for which there exists $\gamma_2\in (-\frac{n+1}{2},-\frac{n}{2})$ such that \[\phi^{\gamma_2}(\lambda x)=\lambda^{\gamma_2}\phi^{\gamma_2}(x)\] holds for every $x\in\mathbb{R}^n_+$ and every $\lambda>0$, and such that \[\int_{\mathbb{R}^n_+}|\phi^{\gamma_2}(x)\phi^{\gamma_2}(1+x)|\d{x}<\infty.\] The corresponding process $G$ is a well-defined $\gamma$-self-similar process with stationary increments (here $\gamma=\gamma_1+\gamma_2+\frac{n}{2}+1$) provided that $-1<-\gamma_2-\frac{n}{2}-1<\gamma_1<-\gamma_2-\frac{n}{2}<\frac{1}{2}$; see \cite[Theorem~3.27]{BaiTaq14}. This family of stochastic processes generalizes the one treated above; however, without assuming the product structure in \eqref{eq:A_Hermite_general} as in \eqref{eq:A_Hermite}, it remains unclear whether our results can be applied in this case.
\end{remark}

\section{Tensor calculus: Cancellations of tensors and expansion formula}
\label{sec:tensor_calc}

In this section, we review elements of tensor calculus and an explicit formula for Wiener chaos expansion of products of random variables. We adopt the following

\begin{convention}
When working in the field of integers, $[a,b]$ shall denote the set $\{i\in\mathbb Z:\,a\le i \mbox{ and } i\le b\}$. When working in the filed of reals, $[a,b]$ shall denote the set $\{t\in\mathbb R:\,a\le t\mbox{ and }t\le b\}$.
\end{convention}

\subsection{Cancellation of tensors}

\subsubsection{Motivating example}
\label{sec:ME1}

Consider a product of Hilbert spaces 
\[
H^{\otimes_2^5}\times H^{\otimes_2^5}\times H^{\otimes_2^3}\times H^{\otimes_2^1}\times H^{\otimes_2^4},
\]
define sets of indices corresponding to the coordinates of this space
\[
I_1=\{1,2,3,4,5\},\quad I_2=\{6,7,8,9,10\},\quad I_3=\{11,12,13\},\quad I_4=\{14\},\quad I_5=\{15,16,17,18\},
\]
and consider a set $V$ of unordered pairs of indices that obey the following rules:

\begin{enumerate}
\item Every index is in one pair at most.
\item No index is paired with an index from the same set nor with itself.
\end{enumerate}

\begin{figure}[H]
\centering
\begin{tikzpicture}[
dot/.style = {circle, fill, minimum size=#1, inner sep=0pt, outer sep=0pt},
dot/.default = 3pt 
]  
\node[draw,circle,minimum size=70,label=above left:$I_1$] (CircleNode) at (-2.5,1.5) {};
\node[dot,label=below:{\scriptsize $1$}] at (-2.5,0.7) {};
\node[dot,label=right:{\scriptsize $2$}] at (-3.2,1) {};
\node[dot,label=below:{\scriptsize $3$},name=3] at (-3,2.1) {};
\node[dot,label=below:{\scriptsize $4$}] at (-2,1.9) {};
\node[dot,label=below:{\scriptsize $5$},name=5] at (-1.8,1.3) {};

\node[draw,circle,minimum size=60,label=above right:$I_2$] (CircleNode) at (0,3) {};
\node[dot,label=below:{\scriptsize $6$},name=6] at (-0.6,3.3) {};
\node[dot,label=right:{\scriptsize $7$}] at (0,3.5) {};
\node[dot,label=right:{\scriptsize $8$}] at (0.6,3) {};
\node[dot,label=above:{\scriptsize $9$}, name=9] at (0.5,2.4) {};
\node[dot,label=above:{\scriptsize $10$},name=10] at (0,2.4) {};

\node[draw,circle,minimum size=40,label=right:$I_3$] (CircleNode) at (2.5,2) {};
\node[dot,label=below:{\scriptsize $11$},name=11] at (2.2,2.3) {};
\node[dot,label=below:{\scriptsize $12$}] at (2.9,2.3) {};
\node[dot,label=right:{\scriptsize $13$},name=13] at (2.5,1.6) {};

\node[draw,circle,minimum size=30,label=below right:$I_4$] (CircleNode) at (2,0) {};
\node[dot,label=below:{\scriptsize $14$},name=14] at (2,0) {};

\node[draw,circle,minimum size=50,label=below left:$I_5$] (CircleNode) at (-0.5,-0.5) {};
\node[dot,label=below:{\scriptsize $15$}] at (0,-0.7) {};
\node[dot,label=below:{\scriptsize $16$}] at (-0.6,-0.9) {};
\node[dot,label=below:{\scriptsize $17$}] at (-0.9,-0.1) {};
\node[dot,label=below:{\scriptsize $18$}] at (0,0) {};

\draw[red] (3)--(6);
\draw[blue] (5)--(10);
\draw[green] (9)--(11);
\draw[brown] (13)--(14);
\end{tikzpicture}
\caption{Example of tensor cancellation.}
\label{fig:1}
\end{figure}
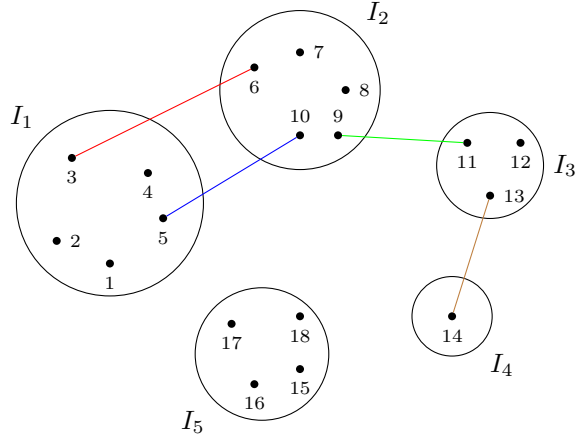

For example, one can consider 
	\[V=\{\{{\color{red}3},{\color{red}6}\},\{{\color{blue} 5},{\color{blue}10}\},\{{\color{green} 9},{\color{green}11}\}, \{{\color{brown}13},{\color{brown} 14}\}\}.\] The example can be graphically visualized as in \autoref{fig:1}.

For the element
\[
(h_{1}\otimes h_{2}\otimes h_{3}\otimes h_{4}\otimes h_{5})\times(h_{6}\otimes h_{7}\otimes h_{8}\otimes h_{9}\otimes h_{10})\times(h_{11}\otimes h_{12}\otimes h_{13})\times h_{14}\times(h_{15}\otimes h_{16}\otimes h_{17}\otimes h_{18}),
\]
the ($V$-)cancellation is done as follows: The vectors with indices forming a pair in $V$ are multiplied and the remaining vectors are shrunken, i.e.

\begin{figure}[H]
\centering
\begin{tikzpicture}
\draw[red, thick] (1.85,0.45) -- (1.85,0.55) -- (4,0.55) -- (4,0.45);
\draw[blue, thick] (3.3,0.45) -- (3.3,0.65) -- (7,0.65) -- (7,0.45);
\draw[green, thick] (6.2,0.45) -- (6.2,0.55) -- (7.85,0.55) -- (7.85,0.45);
\draw[brown, thick] (9.5,0.45) -- (9.5,0.55) -- (10.45,0.55) -- (10.45,0.45);
\node[label={[label distance=0]right:\small $\underbrace{h_{1}\otimes h_{2}\otimes {\color{red}h_{3}}\otimes h_{4}\otimes {\color{blue}h_{5}}}_{I_1}\times\underbrace{{\color{red}h_{6}}\otimes h_{7}\otimes h_{8}\otimes {\color{green}h_{9}}\otimes {\color{blue}h_{10}}}_{I_2}\times\underbrace{{\color{green}h_{11}}\otimes h_{12}\otimes {\color{brown}h_{13}}}_{I_3}\times\underbrace{\color{brown}h_{14}}_{I_4}\times\underbrace{h_{15}\otimes h_{16}\otimes h_{17}\otimes h_{18}}_{I_5}$}] at (0,0) {};
\end{tikzpicture}
\end{figure}
which results in
\[
\langle{\color{red}h_{3}},{\color{red}h_6}\rangle_{H}\langle{\color{blue} h_{5}},{\color{blue}h_{10}}\rangle_{H}\langle{\color{green} h_{9}},{\color{green}h_{11}}\rangle_{H}\langle{\color{brown} h_{13}},{\color{brown}h_{14}}\rangle_{H} h_{1}\otimes h_{2}\otimes h_{4}\otimes h_{7}\otimes h_{8}\otimes h_{12}\otimes h_{15}\otimes h_{16}\otimes h_{17}\otimes h_{18}.
\]

The cancellation extends to a $5$-linear operator from $H^{\otimes^5_{2}}\times H^{\otimes^5_{2}}\times H^{\otimes^3_{2}}\times H^{\otimes^1_{2}}\times H^{\otimes^4_{2}}$ to $H^{\otimes^{10}_{2}}$.

\subsubsection{General case}
\label{GC1}

Consider $\ell\ge 2$, positive integers $d_1,\dots,d_\ell$ and decompose 
	\[
		\{1,\dots,N\}
	\]
to subsequent intervals $I_1,\dots,I_\ell$ of lengths $d_1,\dots,d_\ell$, respectively, where $N=d_1+\dots+d_\ell$. 

\begin{definition}\label{ap}
A set $V$ of unordered pairs of numbers in $\{1,\dots,N\}$ is said to be {\it a set of admissible pairs} for intervals $I_1,\dots,I_{\ell}$ if it is either empty, or it can be enumerated as
	\begin{equation}
	\label{enumV}
			V=\{\{m_1,n_1\},\dots,\{m_k,n_k\}\}
	\end{equation}
for some $k\ge 1$ where 
	\begin{enumerate}[label=(F\arabic*)]
		\item\label{ass:F1} $m_1,\dots,m_k,n_1,\dots,n_k$ are all distinct, i.e. $|\{m_1,\dots,m_k,n_1,\dots,n_k\}|=2k$, and
		\item\label{ass:F2} $m_j$ and $n_j$ do not belong to the same interval $I_1,\dots,I_\ell$ for every $j\in\{1,\dots,k\}$.
	\end{enumerate}
\end{definition}

For a set $V$ of admissible pairs for intervals $I_1,\ldots, I_{\ell}$ with $|V|=k$, we define 
$V_*=\{1,\dots,N\}$ if $V=\emptyset$ or, if $V\ne\emptyset$,
	\begin{equation}
	\label{vstar}
		V_*=\{1,\dots,N\}\setminus\{m_1,\dots,m_k,n_1,\dots,n_k\}.
	\end{equation}
If $2k<N$, we enumerate $V^*$ in an increasing order as 
	\begin{equation}
	\label{def_o}
		\{o_1,\dots,o_{N-2k}\}.
	\end{equation}
The cancellation operator is then defined for $h_1,\ldots, h_N\in H$ as follows: If $V=\emptyset$, then
	\[
		R^{V}(\otimes_{i_1\in I_1}h_{i_1},\dots,\otimes_{i_\ell\in I_\ell}h_{i_\ell})=h_{o_1}\otimes\dots\otimes h_{o_N}=h_1\otimes\dots\otimes h_N,
	\]
if $0<2k<N$, then
	\[
		R^{V}(\otimes_{i_1\in I_1}h_{i_1},\dots,\otimes_{i_\ell\in I_\ell}h_{i_\ell})=\prod_{j=1}^k\langle h_{m_j},h_{n_j}\rangle_{H} h_{o_1}\otimes\dots\otimes h_{o_{N-2k}},
	\]
and if $2k=N$, then
	\[
		R^{V}(\otimes_{i_1\in I_1}h_{i_1},\dots,\otimes_{i_\ell\in I_\ell}h_{i_\ell})=\prod_{j=1}^k\langle h_{m_j},h_{n_j}\rangle_{H}.
	\]
Due to symmetry, the definition of $R^{V}$ is independent of the enumeration of the set $V$ and $R^{V}$ extends to a unique $\ell$-linear operator
	\[
		R^{V}:H^{\otimes^{d_1}_{2}}\times\dots\times H^{\otimes^{d_\ell}_{2}}\to H^{\otimes^{N-2|V|}_{2}}.
	\]

\subsubsection{Permutations}
\label{sec_perm}

It is often convenient to use a relation between cancellations and permutations. Under the assumptions made in \autoref{GC1}, let us assume additionally that $\pi_j$ is a permutation on $\{1,\dots,d_j\}$ for every $j\in\{1,\dots,\ell\}$ and define
	\[
		\pi(i)=s_j+\pi_j(i-s_j),\qquad i\in I_j,\quad j\in\{1,\dots,\ell\},
	\]
where $I_j=s_j+\{1,\dots,d_j\}$ for $j\in\{1,\dots,\ell\}$. Then $\pi$ is a permutation on $\{1,\dots,N\}$ such that $\pi|_{I_j}$ is a permutation on $I_j$ for every $j\in\{1,\dots,\ell\}$. We define $V^\pi=\emptyset$ if $V=\emptyset$ and
	\[
		V^\pi=\{\{\pi(m_1),\pi(n_1)\},\dots,\{\pi(m_k),\pi(n_k)\}\}
	\]
otherwise, with the notation from \eqref{enumV}. It is clear that $V^\pi$ is admissible for the intervals $I_1,\dots,I_{\ell}$. If $2k<N$, let $\sigma$ be the permutation on $\{1,\dots,N-2k\}$ such that $\pi\circ o\circ\sigma$ is increasing. Then
	\begin{equation}
	\label{prem_prop}
		R^{V^\pi}(A_1,\dots,A_\ell)=P_\sigma R^V(P_{\pi_1}A_1,\dots,P_{\pi_\ell}A_\ell), \qquad A_1\in H^{\otimes_2^{d_1}},\ldots, A_\ell \in H^{\otimes_2^{d_\ell}},
	\end{equation}
where the permutation operator $P_\theta:H^{\otimes^n_2}\to H^{\otimes^n_2}$ is defined in a standard manner by
	\begin{equation}
	\label{prem_def}
		P_\theta(h_1\otimes\dots\otimes h_n)=h_{\theta_1}\otimes\dots\otimes h_{\theta_n}\quad\text{if}\quad n\in\mathbb N\quad\text{and}\quad P_\emptyset(t)=t\qquad\text{if}\quad n=0.
	\end{equation}
If $2k=N$ then \eqref{prem_prop} still holds with $\sigma=\emptyset$ and $P_\emptyset=1$.

\subsection{Composition of cancellations}
\subsubsection{Motivating example - continued}

\begin{figure}[H]
\centering
\begin{tikzpicture}[
dot/.style = {circle, fill, minimum size=#1, inner sep=0pt, outer sep=0pt},
dot/.default = 3pt 
]  
\node[dot,label=below:{\scriptsize $1$}] at (-2.5,0.7) {};
\node[dot,label=right:{\scriptsize $2$}, name=2] at (-3.2,1) {};
\node[dot,name=3] at (-3,2.1) {};
\node[dot,label=below:{\scriptsize $3$}] at (-2,1.9) {};
\node[dot,name=5] at (-1.8,1.3) {};

\node[dot,name=6] at (-0.6,3.3) {};
\node[dot,label=right:{\scriptsize $4$}] at (0,3.5) {};
\node[dot,label=right:{\scriptsize $5$}] at (0.6,3) {};
\node[dot,name=9] at (0.5,2.4) {};
\node[dot,name=10] at (0,2.4) {};

\node[dot,name=11] at (2.2,2.3) {};
\node[dot,label=below:{\scriptsize $6$}] at (2.9,2.3) {};
\node[dot,name=13] at (2.5,1.6) {};

\node[dot,name=14] at (2,0) {};

\node[dot,label=below:{\scriptsize $7$}] at (0,-0.7) {};
\node[dot,label=below:{\scriptsize $8$}] at (-0.6,-0.9) {};
\node[dot,label=below:{\scriptsize $9$},name=17] at (-0.9,-0.1) {};
\node[dot,label=below:{\scriptsize $10$}] at (0,0) {};

\node[draw,cyan,circle,minimum size=40,label=below left:$I_6$] (CircleNode) at (-4,-0.5) {};
\node[dot,label=below left:{\scriptsize $11$},name=19] at (-4.1,-0.2) {};
\node[dot,label=below:{\scriptsize $12$}] at (-4.1,-0.7) {};
\node[dot,label=below:{\scriptsize $13$}, name=21] at (-3.7,-0.4) {};

\node[draw,circle,minimum size=200] (CircleNode) at (0,1.5) {};

\draw[dashed, red] (3)--(6);
\draw[dashed, red] (5)--(10);
\draw[dashed, red] (9)--(11);
\draw[dashed, red] (13)--(14);
\draw[cyan] (19)--(2);
\draw[cyan] (21)--(17);
\end{tikzpicture}
\caption{Example of tensor cancellation with composition - new enumeration.}
\label{fig:2}
\end{figure}
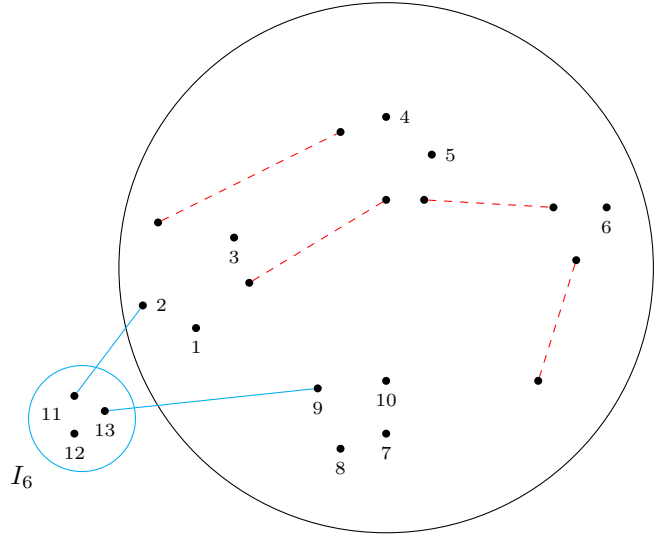

Let us start again by a sequel of the previous example in \autoref{sec:ME1} after cancellation. We renumerate the remaining indices in the ascending order and we add a new interval $I_6=\{11,12,13\}$ with two more pairs $V^\prime=\{\{2,11\},\{9,13\}\}$ obeying rules \ref{ass:F1} and \ref{ass:F2}. This situation is now depicted in \autoref{fig:2}.

\begin{figure}[H]
\centering
\begin{tikzpicture}[
dot/.style = {circle, fill, minimum size=#1, inner sep=0pt, outer sep=0pt},
dot/.default = 3pt 
]  
\node[draw,circle,minimum size=70,label=above left:$I_1$] (CircleNode) at (-2.5,1.5) {};
\node[dot,label=below:{\scriptsize $1$}] at (-2.5,0.7) {};
\node[dot,label=right:{\scriptsize $2$}, name=2] at (-3.2,1) {};
\node[dot,label=below:{\scriptsize $3$},name=3] at (-3,2.1) {};
\node[dot,label=below:{\scriptsize $4$}] at (-2,1.9) {};
\node[dot,label=below:{\scriptsize $5$},name=5] at (-1.8,1.3) {};

\node[draw,circle,minimum size=60,label=above right:$I_2$] (CircleNode) at (0,3) {};
\node[dot,label=below:{\scriptsize $6$},name=6] at (-0.6,3.3) {};
\node[dot,label=right:{\scriptsize $7$}] at (0,3.5) {};
\node[dot,label=right:{\scriptsize $8$}] at (0.6,3) {};
\node[dot,label=above:{\scriptsize $9$}, name=9] at (0.5,2.4) {};
\node[dot,label=above:{\scriptsize $10$},name=10] at (0,2.4) {};

\node[draw,circle,minimum size=40,label=right:$I_3$] (CircleNode) at (2.5,2) {};
\node[dot,label=below:{\scriptsize $11$},name=11] at (2.2,2.3) {};
\node[dot,label=below:{\scriptsize $12$}] at (2.9,2.3) {};
\node[dot,label=right:{\scriptsize $13$},name=13] at (2.5,1.6) {};

\node[draw,circle,minimum size=30,label=below right:$I_4$] (CircleNode) at (2,0) {};
\node[dot,label=below:{\scriptsize $14$},name=14] at (2,0) {};

\node[draw,circle,minimum size=50,label=below left:$I_5$] (CircleNode) at (-0.5,-0.5) {};
\node[dot,label=below:{\scriptsize $15$}] at (0,-0.7) {};
\node[dot,label=below:{\scriptsize $16$}] at (-0.6,-0.9) {};
\node[dot,label=below:{\scriptsize $17$},name=17] at (-0.9,-0.1) {};
\node[dot,label=below:{\scriptsize $18$}] at (0,0) {};

\node[draw,cyan,circle,minimum size=40,label=below left:$I_6$] (CircleNode) at (-4,-0.5) {};
\node[dot,label=below left:{\scriptsize $19$},name=19] at (-4.1,-0.2) {};
\node[dot,label=below:{\scriptsize $20$}] at (-4.1,-0.7) {};
\node[dot,label=below:{\scriptsize $21$}, name=21] at (-3.7,-0.4) {};

\draw[red] (3)--(6);
\draw[red] (5)--(10);
\draw[red] (9)--(11);
\draw[red] (13)--(14);
\draw[cyan] (19)--(2);
\draw[cyan] (21)--(17);
\end{tikzpicture}
\caption{Example of tensor cancellation with composition - original enumeration.}
\label{fig:3}
\end{figure}
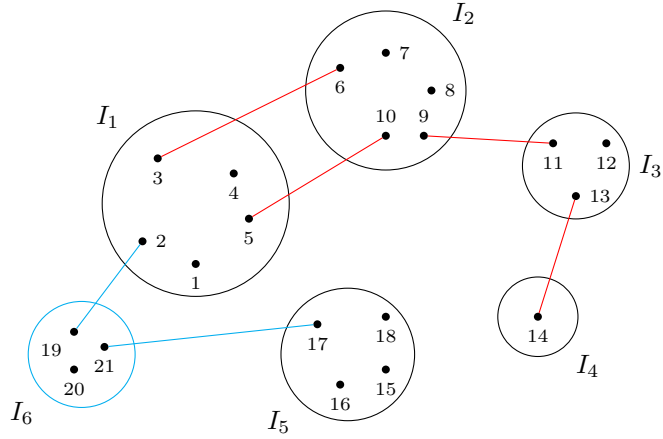

In the original picture (with the original enumeration), the new indices are enumerated as $I_6=\{19,20,21\}$ and the new pairs correspond to the set $V^{\prime\prime}=\{\{2,19\},\{17,21\}\}$ as shown in \autoref{fig:3}.

The extended cancellation $R^{V\cup V''}: H^{\otimes^5_2}\times H^{\otimes^5_2}\times H^{\otimes^3_2}\times H^{\otimes^1_2}\times H^{\otimes^4_2}\times H^{\otimes_2^3} \rightarrow H^{\otimes^{9}_{2}}$ then satisfies the formula
\[
R^{V\cup V^{\prime\prime}}(A_1,A_2,A_3,A_4,A_5,A_6)=R^{V^\prime}(R^V(A_1,A_2,A_3,A_4,A_5),A_6).
\]
as

\begin{figure}[H]
\centering
\begin{tikzpicture}
\draw[red, thick] (1.6,0.4) -- (1.6,0.5) -- (3.35,0.5) -- (3.35,0.4);
\draw[red, thick] (2.8,0.4) -- (2.8,0.6) -- (5.85,0.6) -- (5.85,0.4);
\draw[red, thick] (5.25,0.4) -- (5.25,0.5) -- (6.6,0.5) -- (6.6,0.4);
\draw[red, thick] (8.05,0.4) -- (8.05,0.5) -- (8.8,0.5) -- (8.8,0.4);
\draw[cyan, thick] (1,0.4) -- (1,0.7) -- (12.45,0.7) -- (12.45,0.4);
\draw[cyan, thick] (10.95,0.4) -- (10.95,0.5) -- (13.85,0.5) -- (13.85,0.4);
\node[label={[label distance=0]right:\tiny $\underbrace{h_{1}\otimes {\color{cyan}h_{2}}\otimes {\color{red}h_{3}}\otimes h_{4}\otimes {\color{red}h_{5}}}_{I_1}\times\underbrace{{\color{red}h_{6}}\otimes h_{7}\otimes h_{8}\otimes {\color{red}h_{9}}\otimes {\color{red}h_{10}}}_{I_2}\times\underbrace{{\color{red}h_{11}}\otimes h_{12}\otimes {\color{red}h_{13}}}_{I_3}\times\underbrace{\color{red}h_{14}}_{I_4}\times\underbrace{h_{15}\otimes h_{16}\otimes {\color{cyan}h_{17}}\otimes h_{18}}_{I_5}\times\underbrace{{\color{cyan}h_{19}}\otimes h_{20}\otimes {\color{cyan}h_{21}}}_{I_6}$}] at (0,0) {};
\end{tikzpicture}
\end{figure}

results in
\begin{multline*}
\langle{\color{red}h_{3}},{\color{red}h_6}\rangle_{H}\langle{\color{red} h_{5}},{\color{red}h_{10}}\rangle_{H}\langle{\color{red} h_{9}},{\color{red}h_{11}}\rangle_{H}\langle{\color{red} h_{13}},{\color{red}h_{14}}\rangle_{H}\langle{\color{cyan} h_{2}},{\color{cyan}h_{19}}\rangle_{H}\langle{\color{cyan} h_{17}},{\color{cyan}h_{21}}\rangle_{H} \\ h_{1}\otimes h_{4}\otimes h_{7}\otimes h_{8}\otimes h_{12}\otimes h_{15}\otimes h_{16}\otimes h_{18}\otimes h_{20}.
\end{multline*}

\subsubsection{General case} Under the assumptions in \autoref{GC1}, assume additionally that $2k<N$, consider two intervals 
\[
	I^\prime_1=\{1,\dots,N-2k\},\quad I^\prime_2=\{N-2k+1,\dots,N-2k+d_{\ell+1}\}, 
\]
and consider the interval $I_{\ell+1}=\{N+1,\dots,N+d_{\ell+1}\}$ which completes $I_1,\dots,I_{\ell}$ to a subsequent decomposition of $\{1,\dots,N^\prime\}$ where $N^\prime=N+d_{\ell+1}=d_1+\dots+d_{\ell+1}$. Consider a set of admissible pairs $V^\prime$ for intervals $I^\prime_1$ and $I^\prime_2$. If $V^\prime=\emptyset$, we set $V^{\prime\prime}=\emptyset$. If $V^\prime$ is non-empty, we enumerate it
\[
	V^\prime=\{\{m^\prime_1,n^\prime_1\},\dots,\{m^\prime_{k^\prime},n^\prime_{k^\prime}\}\}
\]
in such a way that $\{m^\prime_1,\dots,m^\prime_{k^\prime}\}\subseteq I^\prime_1$ and $\{n^\prime_1,\dots,n^\prime_{k^\prime}\}\subseteq I^\prime_2$, and we define
$$
m_{k+1}=o(m^\prime_1),\dots, m_{k+k^\prime}=o(m^\prime_{k^\prime}),n_{k+1}=n^\prime_1+2k,\dots,n_{k+k^\prime}=n^\prime_{k^\prime}+2k
$$
and $V^{\prime\prime}=\{\{m_{k+1},n_{k+1}\},\dots,\{m_{k+k^\prime},n_{k+k^\prime}\}\}$.

\begin{lemma}\label{lem_2} The set of pairs $V\cup V^{\prime\prime}$ is admissible for the intervals $I_1,\dots,I_{\ell+1}$ and
\[
R^{V\cup V^{\prime\prime}}(A_1,\dots,A_{\ell+1})=R^{V^\prime}(R^V(A_1,\dots,A_\ell),A_{\ell+1})
\]
holds for every $A_1\in H^{\otimes_2^{d_1}}, \ldots, A_{\ell+1}\in H^{\otimes_2^{d_{\ell+1}}}$.
\end{lemma}

\begin{proof}
If $V^\prime=\emptyset$ then the assertion is trivial. If $V^\prime\ne\emptyset$ but $V=\emptyset$ then $I^\prime_1=\{1,\dots,N\}$, $I^\prime_2=I_{\ell+1}$, $V^\prime=V^{\prime\prime}$ so the assertion is obvious. Now assume that $V\ne\emptyset$ and $V^\prime\ne\emptyset$. The admissibility of the pairs $V\cup V^{\prime\prime}$ is rather straightforward since it is a mere renumeration and, as for the identity, it suffices to show that it holds for elementary tensors, i.e. that
	\begin{multline}
	\label{lhs_1}
		R^{V^\prime}(h_{o_1}\otimes\dots\otimes h_{o_{N-2k}},h_{N+1}\otimes\dots\otimes h_{N+d_{\ell+1}})\\
=\langle h_{m_{k+1}},h_{n_{k+1}}\rangle\dots\langle h_{m_{k+k^\prime}},h_{n_{k+k^\prime}}\rangle h_{o^{\prime\prime}_1}\otimes\dots\otimes h_{o^{\prime\prime}_{N^\prime-2k-2k^\prime}}
	\end{multline}
where $o^{\prime\prime}$ is the increasing enumeration of $\{1,\dots,N^\prime\}\setminus\{m_1,\dots,m_{k+k^\prime},n_1,\dots,n_{k+k^\prime}\}$. If we define
	\[
		h^\prime_1=h_{o_1},\dots,h^\prime_{N-2k}=h_{o_{N-2k}},h^\prime_{N-2k+1}=h_{N+1},\dots,h^\prime_{N-2k+d_{\ell+1}}=h_{N+d_{\ell+1}},
	\]
then the left-hand-side of \eqref{lhs_1} is equal to 
	\begin{equation}
	\label{lhs_2}
		\langle h^\prime_{m^\prime_1},h^\prime_{n^\prime_1}\rangle\dots\langle h^\prime_{m^\prime_{k^\prime}},h^\prime_{n^\prime_{k^\prime}}\rangle h^\prime_{o^\prime_1}\otimes\dots\otimes h^\prime_{o^\prime_{N^\prime-2k-2k^\prime}}
	\end{equation}
where $o^\prime$ is the increasing enumeration of $\{1,\dots,N^\prime-2k\}\setminus\{m^\prime_1,\dots,m^\prime_{k^\prime},n^\prime_1,\dots,n^\prime_{k^\prime}\}$. Since $m^\prime_j\le N-2k<n^\prime_j$, we have $h^\prime_{m^\prime_j}=h_{o(m^\prime_j)}=h_{m_{k+j}}$ and $h^\prime_{n^\prime_j}=h_{n^\prime_j+2k}=h_{n_{k+j}}$. Let us also define
	\[
		L=(o_1,\dots,o_{N-2k},N+1,\dots,N+d_{\ell+1}).
	\]
Then 
	\[
		L:\{1,\dots,N^\prime-2k\}\to\{1,\dots,N^\prime\}\setminus\{m_1,\dots,m_k,n_1,\dots n_k\}
	\]
is an increasing bijection and it can be checked that 
	\[
		\{m_{k+1},\dots,m_{k+k^\prime},n_{k+1},\dots,n_{k+k^\prime}\}\cap\operatorname{Rng}\,(L(o^\prime))=\emptyset.
	\]
Thus it is seen that
	\[
L(o^\prime):\{1,\dots,N^\prime-2k-2k^\prime\}\to\{1,\dots,N^\prime\}\setminus\{m_1,\dots,m_{k+k^\prime},n_1,\dots n_{k+k^\prime}\}
	\]
is an increasing bijection. But such bijection is exactly one, hence $o^{\prime\prime}=L(o^\prime)$. Now, $h^\prime_i=h_{L(i)}$ for $i\in\{1,\dots,N^\prime-2k\}$ so $h^\prime_{o^\prime_j}=h_{L(o^\prime_j)}=h_{o^{\prime\prime}_j}$ for every $j\in\{1,\dots,N^\prime-2k-2k^\prime\}$ which proves that \eqref{lhs_2} coincides with the right-hand side of \eqref{lhs_1}.
\end{proof}

We will prove the following two results using the composition \autoref{lem_2}.

\begin{lemma}
\label{mult_est} 
If $A_1\in H^{\otimes^{d_1}},\dots, A_\ell\in H^{\otimes^{d_\ell}}$, then
	\[
		\|R^{V}(A_1,\dots,A_\ell)\|_{H^{\otimes_2^{N-2|V|}}}\le\|A_1\|_{H^{\otimes_2^{d_1}}}\dots\|A_\ell\|_{H^{\otimes_2^{d_\ell}}}.
	\]
\end{lemma}

\begin{proof}
Let us proceed by induction on $\ell\ge 2$. Let $U$ be an admissible set of pairs for intervals $I_1,\dots,I_{\ell+1}$ and define
	\[
		V=\{\{m,n\}\in U:\,m\notin I_{\ell+1}\,\,\text{and}\,\,n\notin I_{\ell+1}\},\qquad V^{\prime\prime}=U\setminus V,
	\]
with cardinalities $|V|=k$ 
and $|V^{\prime\prime}|=k^\prime$, 
respectively. We distinguish two cases. First, if $V^{\prime\prime}=\emptyset$, then define $V^\prime=\emptyset$. If $V^{\prime\prime}\neq\emptyset$, then necessarily $2k<N$ where $N=d_1+\dots+d_\ell$ and $N^\prime=N+d_{\ell+1}$, so we can enumerate 
	\[
		V^{\prime\prime}=\{\{m_{k+1},n_{k+1}\},\dots,\{m_{k+k^\prime},n_{k+k^\prime}\}\}
	\]
such that $m_{k+j}\le N<n_{k+j}$ for every $j\in[1,k^\prime]$. Since $U$ is admissible, $m_{k+1},\dots,m_{k+k^\prime}$ belong to $V_*$, hence we can define $m^\prime_j=o^{-1}(m_{k+j})\in[1,N-2k]=I^\prime_1$ 
and $n^\prime_j=n_{k+j}-2k\in [N-2k+1,N-2k+d_{\ell+1}]=I^\prime_2$ for $j\in [1,k^\prime]$, and $V^\prime=\{\{m_j^\prime,n_j^\prime\},j\in [1,k^\prime]\}$. After this construction, \autoref{lem_2} is applicable. If $V^{\prime\prime}=\emptyset$, then $U=V$ and \[R^U(A_1,\dots,A_{\ell+1})=R^V(A_1,\dots,A_\ell)\otimes A_{\ell+1}\] holds if $2k<N$ or \[R^U(A_1,\dots,A_{\ell+1})=R^V(A_1,\dots,A_\ell)A_{\ell+1}\] holds if $2k=N$. If $V^{\prime\prime}\ne\emptyset$, then \[R^U(A_1,\dots,A_{\ell+1})=R^{V^\prime}(R^V(A_1,\dots,A_\ell),A_{\ell+1}),\] so the claim follows by the induction step. Let therefore $\ell=2$. Since $R^\emptyset(A_1,A_2)=A_1\otimes A_2$, it suffices to assume that $V\ne\emptyset$. Express
	\[
		A_1=\sum_{i_1,\dots,i_{d_1}}\alpha^1_{i_1,\dots,i_{d_1}}e_{i_1}\otimes\dots\otimes e_{i_{d_1}},\qquad A_2=\sum_{i_{d_1+1},\dots,i_{d_1+d_2}}\alpha^2_{i_{d_1+1},\dots,i_{d_1+d_2}}e_{i_{d_1+1}}\otimes\dots\otimes e_{i_{d_1+d_2}}
	\]
for some orthonormal system $\{e_\beta\}$ in $H$ and let $V=\{\{j,d_1+j\}:\,j\in\{1,\dots,k\}\}$ for some $k\le\min\,\{d_1,d_2\}$. Then, if $k<\min\,\{d_1,d_2\}$, we have that
	\[
		R^V(A_1,A_2)=\sum_{i,j,l}\alpha^1_{i,j}\alpha^2_{i,l}e_{j_1}\otimes\dots\otimes e_{j_{d_1-k}}\otimes e_{l_1}\otimes\dots\otimes e_{l_{d_2-k}}
	\]
so the estimate
	\begin{equation}
	\label{2_est_ord}
		\|R^V(A_1,A_2)\|^2_{H^{\otimes_2^{d_1+d_2-2k}}} =\sum_{j,l}\left(\sum_i\alpha^1_{i,j}\alpha^2_{i,l}\right)^2\le\|A_1\|^2_{H^{\otimes_2^{d_1}}}\|A_2\|^2_{H^{\otimes_2^{d_2}}}
	\end{equation}
follows by the Cauchy-Schwarz inequality. If $k=\min\,\{d_1,d_2\}$, the same estimation is obtained analogously. In the general case, the result is obtaind from \eqref{prem_prop} and \eqref{2_est_ord} by reordering the set $V$ via a suitable permutation.
\end{proof}

\begin{corollary} The operator $R^{V}$ extends to a continuous $\ell$-linear operator
	\[
		R^{V}:H^{\otimes_2^{d_1}}\times\dots\times H^{\otimes_2^{d_\ell}}\to H^{\otimes_2^{N-2|V|}}.
	\]
\end{corollary}

Now we are going to study the behaviour of the operator $R$ further. Recall the numbers $s_1,\dots,s_\ell$ from \autoref{sec_perm} that were defined in such a way that
	\[
		I_1=s_1+\{1,\dots,d_1\},\dots,I_\ell=s_\ell+\{1,\dots,d_\ell\},
	\]
holds, consider the set $V_*$ from \eqref{vstar}, and define
	\[
		V_*^{(1)}=(V_*\cap I_1)-s_1,\dots V_*^{(\ell)}=(V_*\cap I_\ell)-s_\ell.
	\]
In this way, $V_*^{(j)}\subseteq\{1,\dots,d_j\}$ for every $j\in\{1,\dots,\ell\}$ and these sets are actually the traces of $V_*$ on $I_j$, but renumbered such that each interval begins with $1$. Finally, define
	\[
		V_j=\{\{i,d_j+i\}:\,i\in V_*^{(j)}\},\qquad j\in\{1,\dots,\ell\}.
	\]

\begin{lemma}\label{mult_est_2} If $A_1,B_1\in H^{\otimes_2^{d_1}},\dots A_\ell,B_\ell\in H^{\otimes_2^{d_\ell}}$ then
	\begin{multline*}
	|\langle R^V(A_1,\dots,A_\ell),R^{V}(B_1,\dots,B_\ell)\rangle_{H^{\otimes_2^{N-2|V|}}}|\\ \le\|R^{V_1}(A_1,B_1)\|_{H^{\otimes_2^{N-2|V_1|}}}\dots\|R^{V_\ell}(A_\ell,B_\ell)\|_{H^{\otimes_2^{N-2|V_\ell|}}}.
	\end{multline*}
\end{lemma}

\begin{proof}
If $V=\emptyset$, then
	\begin{align*}
		\langle R^V(A_1,\dots,A_\ell),R^{V}(B_1,\dots,B_\ell)\rangle _{H^{\otimes_2^{N-2|V|}}}
		&=\langle A_1\otimes\dots\otimes A_\ell,B_1\otimes\dots\otimes B_\ell\rangle_{H^{\otimes_2^{N-2|V|}}}\\
		&=\langle A_1,B_1\rangle_{H^{\otimes_2^{d_1}}}\dots\langle A_\ell,B_\ell\rangle_{H^{\otimes_2^{d_\ell}}}\\
		&=R^{V_1}(A_1,B_1)\dots R^{V_\ell}(A_\ell,B_\ell).
	\end{align*}
Assume therefore that $V\ne\emptyset$. First assume additionally that, for every $j\in\{1,\dots,\ell\}$, the sets $V^{(j)}_*$ are either empty or $V^{(j)}_*=\{\lambda_j+1,\dots,d_j\}$ for some $\lambda_j$. We can cover both cases simultaneously by admitting that $\lambda_j\in\{0,\dots,d_j\}$.
Now, consider an expansion
	\[
		A_j=\sum_{i_1,\dots,i_{d_j}}a^j_{i_1,\dots,i_{d_j}}e_{i_1}\otimes\dots\otimes e_{i_{d_j}},\qquad B_j=\sum_{i_1,\dots,i_{d_j}}b^j_{i_1,\dots,i_{d_j}}e_{i_1}\otimes\dots\otimes e_{i_{d_j}}
	\]
for some orthonormal system $\{e_\gamma\}$ in $H$. Then 
\begin{align*}
R^V(A_1,\dots,A_\ell)&=\sum_i\left(\prod_{\alpha=1}^\ell a^\alpha_{i_{s_\alpha+1},\dots,i_{s_\alpha+d_\alpha}}\right)\left(\prod_{\beta=1}^k\delta_{i_{m_\beta},i_{n_\beta}}\right)\bigotimes_{\gamma=1}^\ell\left(e_{i_{s_\gamma+\lambda_\gamma+1}}\otimes\dots\otimes e_{i_{s_\gamma+d_\gamma}}\right)
\\
R^V(B_1,\dots,B_\ell)&=\sum_i\left(\prod_{\alpha=1}^\ell b^\alpha_{i_{s_\alpha+1},\dots,i_{s_\alpha+d_\alpha}}\right)\left(\prod_{\beta=1}^k\delta_{i_{m_\beta},i_{n_\beta}}\right)\bigotimes_{\gamma=1}^\ell\left(e_{i_{s_\gamma+\lambda_\gamma+1}}\otimes\dots\otimes e_{i_{s_\gamma+d_\gamma}}\right).
\end{align*}
Define, for $\alpha\in\{1,\dots,\ell\}$,
	\[
		\varrho^\alpha(i_1,\dots,i_{\lambda_\alpha},j_1,\dots,j_{\lambda_\alpha})=\sum_{u_{\lambda_\alpha+1},\dots,u_{d_\alpha}}a^\alpha_{i_1,\dots,i_{\lambda_\alpha},u_{\lambda_\alpha+1},\dots,u_{d_\alpha}}b^\alpha_{j_1,\dots,j_{\lambda_\alpha},u_{\lambda_\alpha+1},\dots,u_{d_\alpha}}
	\]
if $\lambda_\alpha\in\{1,d_\alpha-1\}$ (and if $d_\alpha>1$) and with obvious modifications if $\lambda_\alpha\in\{0,d_\alpha\}$. Then, with the convention that a product over the empty set equals $1$, we have
	\begin{align*}
		\MoveEqLeft[3]\langle R^V(A_1,\dots,A_\ell),R^V(B_1,\dots,B_\ell)\rangle \\
			 ={} & \sum_{i}\sum_j\left(\prod_{\alpha=1}^\ell a^\alpha_{i_{s_\alpha+1},\dots,i_{s_\alpha+d_\alpha}}b^\alpha_{j_{s_\alpha+1},\dots,j_{s_\alpha+d_\alpha}}\prod_{c=\lambda_\alpha+1}^{d_\alpha}\delta_{i_{s_\alpha+c},j_{s_\alpha+c}}\right)\left(\prod_{\beta=1}^k\delta_{i_{m_\beta},i_{n_\beta}}\delta_{j_{m_\beta},j_{n_\beta}}\right)
\\
			={} & \sum_{\substack{i_{s_1+1},\dots,i_{s_1+\lambda_1}\\\dots\\i_{s_\ell+1},\dots,i_{s_\ell+\lambda_\ell}}}\sum_{\substack{j_{s_1+1},\dots,j_{s_1+\lambda_1}\\\dots\\j_{s_\ell+1},\dots,j_{s_\ell+\lambda_\ell}}}\left(\prod_{\alpha=1}^\ell \varrho^\alpha_{i_{s_\alpha+1},\dots,i_{s_\alpha+\lambda_\alpha},j_{s_\alpha+1},\dots,j_{s_\alpha+\lambda_\alpha}}\right)\left(\prod_{\beta=1}^k\delta_{i_{m_\beta},i_{n_\beta}}\delta_{j_{m_\beta},j_{n_\beta}}\right)
\\
			={} &\sum_{\substack{i_{m_1},j_{m_1},\dots,i_{m_k},j_{m_k}\\i_{n_1},j_{n_1},\dots,i_{n_k},j_{n_k}}}\left(\prod_{\alpha=1}^\ell \varrho^\alpha_{i_{s_\alpha+1},\dots,i_{s_\alpha+\lambda_\alpha},j_{s_\alpha+1},\dots,j_{s_\alpha+\lambda_\alpha}}\right)\left(\prod_{\beta=1}^k\delta_{i_{m_\beta},i_{n_\beta}}\delta_{j_{m_\beta},j_{n_\beta}}\right)
	\end{align*}
since
	\[
		\{s_\alpha+l:\,1\le l\le\lambda_\alpha,\,1\le\alpha\le\ell\}=\{m_1,\dots,m_k,n_1,\dots,n_k\}.
	\]
Let $s_{\alpha^1_v}+l^1_v=m_v$, $s_{\alpha^2_v}+l^2_v=n_v$ be the assignment. Then $\alpha^1_v\ne\alpha^2_v$ by \ref{ass:F2} and each of the $4k$ variables $i_{m_1},j_{m_1},\dots,i_{m_k},j_{m_k},i_{n_1},j_{n_1},\dots,i_{n_k},j_{n_k}$ appears only once in the formula above, so we can apply the Cauchy-Schwarz inequality in the form
	\begin{equation}
	\label{casc}
		\sum_{\mu\in M}\sum_{\nu\in M}f_\mu g_\nu\delta_{\mu,\nu}\le\|f\|_{\ell_2(M)}\|g\|_{\ell_2(M)}
	\end{equation}
successively on the variables $(i_{m_1},i_{n_1}),(j_{m_1},j_{n_1}),\dots,(i_{m_k},i_{n_k}),(j_{m_k},j_{n_k})$. Eventually, after the last application of \eqref{casc}, we obtain
	\begin{align*}
		\langle R^V(A_1,\dots,A_\ell),R^V(B_1,\dots,B_\ell)\rangle_{H^{\otimes_2^{N-2|V|}}} & \le\|\varrho^1\|_{\ell_2}\dots\|\varrho^\ell\|_{\ell_2} \\ & =\|R^{V_1}(A_1,B_1)\|_{H^{\otimes_2^{d_1}}}\dots\|R^{V_\ell}(A_\ell,B_\ell)\|_{H^{\otimes_2^{d_\ell}}}.
	\end{align*}
If $V\ne\emptyset$ but is general, we reorganize by a permutation $\pi$ to have the form from the above step and use \eqref{prem_prop} to get
	\begin{align*}
		\MoveEqLeft[3]\langle R^V(A_1,\dots,A_\ell),R^V(B_1,\dots,B_\ell)\rangle_{H^{\otimes_2^{N-2|V|}}} & 
		\\ 
			& =\langle R^{V^\pi}(P_{\pi_1^{-1}}A_1,\dots,P_{\pi_\ell^{-1}}A_\ell),R^{V^\pi}(P_{\pi_1^{-1}}B_1,\dots,P_{\pi_\ell^{-1}}B_\ell)\rangle_{H^{\otimes_2^{N-2|V^\pi|}}}
		\\
			&\le\|R^{V^\pi_1}(P_{\pi_1^{-1}}A_1,P_{\pi_1^{-1}}B_1)\|_{H^{\otimes_2^{N-2|V^{\pi}_1|}}}\dots\|R^{V^\pi_\ell}(P_{\pi_\ell^{-1}}A_\ell,P_{\pi_\ell^{-1}}B_\ell)\|_{H^{\otimes_2^{N-2|V^{\pi}_\ell|}}}
		\\
			&=\|R^{(V_1)^{\pi_1^\prime}}(P_{\pi_1^{-1}}A_1,P_{\pi_1^{-1}}B_1)\|_{H^{\otimes_2^{N-2|(V_1)^{\pi_1^\prime}|}}}\dots\|R^{(V_\ell)^{\pi_\ell^\prime}}(P_{\pi_\ell^{-1}}A_\ell,P_{\pi_\ell^{-1}}B_\ell)\|_{H^{\otimes_2^{N-2|(V_\ell)^{\pi_\ell^\prime}|}}}
		\\
			&=\|R^{V_1}(A_1,B_1)\|_{H^{\otimes_2^{N-2|V_1|}}}\dots\|R^{V_\ell}(A_\ell,B_\ell)\|_{H^{\otimes_2^{N-2|V_\ell|}}}.
	\end{align*}
since $V_*^{\pi,(j)}=\pi_j[V_*^{(j)}]$ and $V^\pi_j=(V_j)^{\pi^\prime_j}$ where 
	\[
		\pi^\prime_j=(\pi_j(1),\dots,\pi_j(d_j),d_j+\pi_j(1),\dots,d_j+\pi_j(d_j))
	\]
for every $j\in\{1,\dots,\ell\}$.
\end{proof}

\subsection{Expansion formula}

Let $W$ be an $H$-isonormal Gaussian process and $W_n:H^{\otimes_2^n}\to L^2(\Omega)$, $n\in\mathbb{N}_0$, be the divergence operators as in \autoref{sec:main}. We adopt the following

\begin{notation}
Let $\ell\ge 2$, let $d_1,\dots,d_\ell$ be positive integers, decompose 
	\[
		\{1,\dots,N\}
	\]
to subsequent intervals $I_1,\dots,I_\ell$ of lengths $d_1,\dots,d_\ell$ respectively, where $N=d_1+\dots+d_\ell$ and denote by $\mathcal E_{I_1,\dots,I_{\ell}}$ the system of all admissible sets of pairs $V$ for the the intervals $I_1,\dots,I_{\ell}$ (\autoref{ap}), and $\mathcal E^k_{I_1,\dots,I_{\ell}}$ its subsystem consisting of $V$ such that $|V|=k$.
\end{notation}

\begin{remark} 
We have the following estimate of the cardinality
	\begin{equation}
	\label{card_e}
		|\mathcal E^k_{I_1,\dots,I_{\ell}}|\le\frac{N!}{2^kk!(N-2k)!};
	\end{equation}
see, e.g., \cite[Chapter 1.5, page 16]{Jan97}.
\end{remark}

\begin{theorem}
\label{exp_thm} 
Let $\ell\ge 2$, let $d_1,\dots,d_\ell$ be positive integers, decompose $\{1,\dots,N\}$ to subsequent intervals $I_1,\dots,I_\ell$ of lengths $d_1,\dots,d_\ell$ respectively, where $N=d_1+\dots+d_\ell$. Then
	\[
		W_{d_1}(A_1)\dots W_{d_{\ell}}(A_\ell)=\sum_{V\in\mathcal E_{I_1,\dots,I_\ell}}W_{N-2|V|}(R^V(A_1,\dots,A_\ell))
	\]
holds for every $A_1\in H^{\otimes_2^{d_1}},\dots,A_\ell\in H^{\otimes_2^{d_\ell}}$.
\end{theorem}

\begin{remark}
In \cite[Theorem 7.33]{Jan97}, the theorem is proved in a different language. Below, we provide a proof based on the tensor calculus.
\end{remark}

\begin{proof}
We have that
	\begin{align*}
	\MoveEqLeft[3] W_{d_1}(\underbrace{h_1\otimes\dots\otimes h_{d_1}}_{A_1})W_{d_2}(\underbrace{h_{d_1+1}\otimes\dots\otimes h_{d_1+d_2}}_{A_2}) 
		\\
		&=\sum_{k=0}^{d_1\land d_2}k!\binom{d_1}{k}\binom{d_2}{k} W_{d_1+d_2-2k}(R^{V_k}(h_1\tilde{\otimes}\dots\tilde{\otimes} h_{d_1},h_{d_1+1}\tilde{\otimes}\dots\tilde{\otimes}h_{d_1+d_2}))
		\\
		&=\sum_{\pi\in P_{I_1}}\sum_{\sigma\in P_{I_2}}\sum_{k=0}^{d_1\land d_2}\frac{k!}{d_1!d_2!}\binom{d_1}{k}\binom{d_2}{k} W_{d_1+d_2-2k}(R^{V_k}(h_{\pi_1}\otimes\dots\otimes h_{\pi_{d_1}},h_{\sigma_1}\otimes\dots\otimes h_{\sigma_{d_2}}))
		\\
		&=\sum_{\pi\in P_{I_1}}\sum_{\sigma\in P_{I_2}}\sum_{k=0}^{d_1\land d_2}\frac{k!}{d_1!d_2!}\binom{d_1}{k}\binom{d_2}{k} W_{d_1+d_2-2k}(R^{\{\{\pi_1,\sigma_1\},\dots,\{\pi_k,\sigma_k\}\}}(A_1,A_2))
		\\
		&=\sum_{k=0}^{d_1\land d_2}\sum_{V\in\mathcal E^k_{I_1,I_2}}W_{d_1+d_2-2k}(R^V(A_1,A_2))
		\\
		&=\sum_{V\in\mathcal E_{I_1,I_2}}W_{d_1+d_2-2|V|}(R^V(A_1,A_2))
	\end{align*}
where $V_k=\{(1,d_1+1),\dots,(k,d_1+k)\}$, $P_I$ 
denotes the set of permutations on a set $I$, and $\tilde{\otimes}$ denotes the symmetric tensor product. The first equality follows by, e.g., \cite[Proposition 1.1.3]{Nua06}, the third equality follows from the fact that $W_n(A)=W_n(B)$ for every $A,B\in H^{\otimes^n_2}$ such that $B$ is a permutation of $A$, and the fourth equality follows from the fact that
	\[
		|\{(\pi,\sigma)\in P_{I_1}\times P_{I_2}:\,\{\{\pi_1,\sigma_1\},\dots,\{\pi_k,\sigma_k\}\}=V\}|=k!(d_1-k)!(d_2-k)!
	\]
holds for every $V\in\mathcal E^k_{I_1,I_2}$. We finish the proof by induction on $\ell$ for
	\[
		A_1=h_1\otimes\dots\otimes h_{d_1},\quad A_2=h_{d_1+1}\otimes\dots\otimes h_{d_1+d_2},\quad A_3=h_{d_1+d_2+1}\otimes\dots\otimes h_{d_1+d_2+d_3},\quad \dots
	\]
By induction hypothesis, it holds that
	\[
		W_{d_1}(A_1)\dots W_{d_\ell}(A_\ell)=\sum_{V\in\mathcal E_{I_1,\dots,I_{\ell-1}}}W_{N-d_\ell-2|V|}(R^V(A_1,\dots,A_{\ell-1}))W_{d_\ell}(A_\ell)
	\]
and we must distinguish two cases. If $2|V|<N-d_\ell$, then, by \autoref{lem_2},
	\begin{equation}
	\label{ind_step_1}
		W_{N-d_\ell-2|V|}(R^V(A_1,\dots,A_{\ell-1}))W_{d_\ell}(A_\ell)=\sum_{U\in\mathcal E_{I_1,\dots,I_\ell}(V)}W_{N-2|U|}(R^U(A_1,\dots,A_\ell))
	\end{equation}
where $\mathcal E_{I_1,\dots,I_\ell}(V)$ contains all those elements $U$ of $\mathcal E_{I_1,\dots,I_\ell}$ that complete $V$ by pairs $\{m,n\}$ for which $m$ or $n$ belongs to $I_\ell$. More rigorously, $U$ must satisfy 
\[ U\in \mathcal E_{I_1,\dots,I_\ell}, \quad V\subseteq U, \quad \mbox{and}\quad U\setminus V\in\mathcal E_{I_1\cup\dots\cup I_{\ell-1},I_\ell}.\]

If $2|V|=N-d_\ell$, then \eqref{ind_step_1} holds as well but \autoref{lem_2} is not needed here as $\mathcal E_{I_1,\dots,I_\ell}(V)=\{V\}$. And indeed,
	\begin{align*}
		W_{N-d_\ell-2|V|}(R^V(A_1,\dots,A_{\ell-1}))W_{d_\ell}(A_\ell)&=R^V(A_1,\dots,A_{\ell-1})W_{d_\ell}(A_\ell)
		\\
		&=W_{N-2|V|}(R^V(A_1,\dots,A_\ell))
		\\
		&=\sum_{U\in\mathcal E_{I_1,\dots,I_\ell}(V)}W_{N-2|U|}(R^U(A_1,\dots,A_\ell)).
	\end{align*}
Thus,
	\begin{align*}
		W_{d_1}(A_1)\dots W_{d_\ell}(A_\ell)&=\sum_{V\in\mathcal E_{I_1,\dots,I_{\ell-1}}}\sum_{U\in\mathcal E_{I_1,\dots,I_\ell}(V)}W_{N-2|U|}(R^U(A_1,\dots,A_\ell))
		\\
		&=\sum_{U\in\mathcal E_{I_1,\dots,I_\ell}}W_{N-2|U|}(R^U(A_1,\dots,A_\ell))
	\end{align*}
since $\{\mathcal E_{I_1,\dots,I_\ell}(V):\,V\in\mathcal E_{I_1,\dots,I_{\ell-1}}\}$ is a partition of $\mathcal E_{I_1,\dots,I_\ell}$.
\end{proof}

\section{Proofs}

\subsection{Proof of \autoref{main_thm}}
\label{sec:proof}

Let $\ell\ge 2$ be even and define $N=\ell n$. Then, according to \autoref{exp_thm},
	\[
		[G(s+\delta)-G(s)]^\ell=\sum_{k=0}^{N/2}W_{N-2k}(d^{(k)}_{s,\delta}),\qquad s\in[0,T-\delta],
	\]
where
	\[
		d^{(k)}_{s,\delta}=\sum_{V\in\mathcal E^k_{I_1,\dots,I_\ell}}R^V(A_{s,\delta},\dots,A_{s,\delta}),\qquad s\in[0,T-\delta],
	\]
for $I_j=\{(j-1)n+1,\dots,jn\}$. By \autoref{mult_est}, assumption \ref{ass:G1}, and inequality \eqref{card_e}, there is the estimate
	\begin{equation}
	\label{dsd_est}
		\|d^{(k)}_{s,\delta}\|_{H^{\otimes_2^{N-2k}}}\le\frac{N!\kappa^\ell}{2^kk!(N-2k)!}\delta^{\alpha\ell}
	\end{equation}
from which it follows that the function $s\mapsto d^{(k)}_{s,\delta}$ is integrable and if we define
	\[
		Y_{\ell,\delta}=\|G(\cdot+\delta)-G(\cdot)\|_{L^\ell(0,T-\delta)},\qquad\delta\in (0,T),
	\]
then
	\[
		Y^\ell_{\ell,\delta}=\sum_{k=0}^{N/2}W_{N-2k}(d^{(k)}_{\delta}),\qquad\delta\in (0,T),
	\]
where
	\[
		d^{(k)}_{\delta}=\int_0^{T-\delta}d^{(k)}_{s,\delta}\d{s}=\sum_{V\in\mathcal E^k_{I_1,\dots,I_\ell}}\int_0^{T-\delta}R^V(A_{s,\delta},\dots,A_{s,\delta})\d{s},\qquad\delta\in (0,T).
	\]
Now, for $k=\frac{N}{2}$, we have
	\begin{equation}
	\label{dsd_comp_est}
		|d^{(\frac{N}{2})}_{\delta}|\le\frac{N!}{2^\frac{N}{2}(\frac{N}{2})!}T\kappa^\ell \delta^{\alpha\ell}
	\end{equation}
by estimate \eqref{dsd_est}, and for $k<\frac{N}{2}$, we have
	\begin{align*}
		\MoveEqLeft[1] \|d^{(k)}_{\delta}\|_{H^{\otimes_2^{N-2k}}} &  \\
		& \le\sum_{V\in\mathcal E^k_{I_1,\dots,I_\ell}}\left\|\int_0^{T-\delta}R^V(A_{s,\delta},\dots,A_{s,\delta})\,ds\right\|_{H^{\otimes_2^{N-2k}}}
		\\
		&=\sum_{V\in\mathcal E^k_{I_1,\dots,I_\ell}}\left(\int_0^{T-\delta}\int_0^{T-\delta}\langle R^V(A_{s,\delta},\dots,A_{s,\delta}),R^V(A_{t,\delta},\dots,A_{t,\delta})\rangle_{H^{\otimes_2^{N-2k}}}\d{s}\d{t}\right)^\frac{1}{2}
		\\
		&\le\sum_{V\in\mathcal E^k_{I_1,\dots,I_\ell}}\left(\int_0^{T-\delta}\int_0^{T-\delta}\|R^{V_1}(A_{s,\delta},A_{t,\delta})\|_{H^{\otimes_2^{2n-2|V_1|}}}\dots\|R^{V_\ell}(A_{s,\delta},A_{t,\delta})\|_{H^{\otimes_2^{2n-2|V_\ell|}}}\d{s}\d{t}\right)^\frac{1}{2}
	\end{align*}
by using \autoref{mult_est_2}. Now, we have in fact that $k\leq N-2$ as $N$ is even which means that the set $V_*$ defined by \eqref{vstar} contains at least two elements. It follows that either there are $i,i^\prime\in\{1,\ldots,\ell\}$, $i\neq i^\prime$, such that $V_i\neq \emptyset$ and $V_{i^\prime}\neq\emptyset$ or there is $i\in \{1,\ldots,\ell\}$ such that $V_i$ contains at least two elements. Consequently, the estimate 
	\[ \|R^{V_1}(A_{s,\delta},A_{t,\delta})\|_{H^{\otimes_2^{2n-2|V_1|}}}\dots\|R^{V_\ell}(A_{s,\delta},A_{t,\delta})\|_{H^{\otimes_2^{2n-2|V_\ell|}}} \leq \kappa^{2\ell} \delta^{2\alpha\ell} C_{\delta,\delta}(s,t)\]
is obtained from \ref{ass:G1} by \autoref{mult_est}. This, together with \eqref{card_e}, yields
	\[ \|d_\delta^{(k)}\|_{H^{\otimes_2^{N-2k}}} \leq \frac{N!}{2^kk!(N-2k)!}\kappa^\ell\delta^{\alpha\ell} [F(\delta,\delta)]^\frac{1}{2}.\]
Hence,
	\begin{align*}
		\mathbb E\,(Y^\ell_{\ell,\delta}-\mathbb E\,Y^\ell_{\ell,\delta})^2&=\mathbb E\,\left[\sum_{k=0}^{\frac{N}{2}-1}W_{N-2k}(d^{(k)}_{\delta})\right]^2
		\\
		&=\sum_{k=0}^{\frac{N}{2}-1}\mathbb E\,\left[W_{N-2k}(d^{(k)}_{\delta})\right]^2
		\\
		&\le\sum_{k=0}^{\frac{N}{2}-1}(N-2k)!\|d^{(k)}_{\delta}\|^2_{H^{\otimes_2^{N-2k}}}
		\\
		&\le\kappa^{2\ell}\delta^{2\alpha\ell}N!F(\delta,\delta)\sum_{k=0}^{\frac{N}{2}-1}\frac{(2k)!}{2^{2k}(k!)^2}\binom{N}{2k}
		\\
		&\le \kappa^{2\ell}\delta^{2\alpha\ell}(\ell n)!2^{\ell n} F(\delta,\delta) \numberthis \label{eq:second_moment}
	\end{align*}
since
	\begin{equation}
	\label{fact_est}
		\sup_{k\geq 0}\frac{(2k)!}{2^{2k}(k!)^2}=1.
	\end{equation}
Now, if $q\in (2,\infty)$, then we have 
	\begin{equation}
	\label{eq:hypercontractivity}
	\E|Y_{\ell,\delta}^\ell - \E Y_{\ell,\delta}^\ell|^q \leq N^\frac{q}{2}(q-1)^\frac{qN}{2} [\E(Y_{\ell,\delta}^\ell - \E Y_{\ell,\delta}^\ell)^2]^\frac{q}{2}
	\end{equation}
by the equivalence of moments on a finite Wiener chaos \cite[Corollary 2.8.14]{NouPec12}. Thus the estimate
	\[ \E|Y_{\ell,\delta}^\ell - \E Y_{\ell,\delta}^\ell|^q  
			\leq
		 (\ell n)^\frac{q}{2} (q-1)^\frac{q\ell n}{2} \kappa^{\ell q}\delta^{\alpha\ell q} [(\ell n)!]^\frac{q}{2} 2^{\frac{\ell n q}{2}}[F(\delta,\delta)]^\frac{q}{2}\]
is obtained from \eqref{eq:second_moment}. Let us define $\delta_j=2^{-j}$ in the range $\{j:\,j\ge j_0\}=\{j:\,\delta_j<T\}$ and
	\[
		\bm{C}_q=\left[\sum_{j=j_0}^\infty\sum_{\ell\in 2\mathbb N}\left(\frac{|Y^\ell_{\ell,\delta_j}-\mathbb E\,Y^\ell_{\ell,\delta_j}|}{(2\kappa)^{\ell}(q-1)^\frac{\ell n}{2}\delta_j^{\alpha\ell}\sqrt{2^{\ell n}(\ell n)(\ell n)!}}\right)^{q}\right]^{\frac{1}{q}}.
	\]
Then $[\mathbb E\,\bm C^{q}_q]^\frac{1}{q} \leq [\sum_{j=j_0}^\infty F(\delta_j,\delta_j)]^\frac{1}{2}$ where the sum on the right-hand side is finite by \ref{ass:G3}. Consequently, the constant $\bm C_q$ is almost surely finite and the inequality
	\begin{equation}
	\label{eq:rho}
		|Y^\ell_{\ell,\delta_j}-\mathbb E\,Y^\ell_{\ell,\delta_j}|\le\bm C_q(2\kappa)^{\ell}\delta_j^{\alpha\ell}(q-1)^{\frac{\ell n}{2}}\sqrt{2^{\ell n}(\ell n)(\ell n)!}
	\end{equation}
holds for every $j\ge j_0$ and every even $\ell\geq 2$ almost surely. Moreover, there is the almost sure convergence
	\begin{equation}
	\label{eq:converge_2}
	\lim_{j\rightarrow\infty} \delta_j^{-\alpha\ell}|Y_{\ell,\delta_j}^\ell - \E Y_{\ell,\delta_j}^\ell| =0
	\end{equation}
for every even $\ell\geq 2$. Now
	\[
		\mathbb E\,Y^\ell_{\ell,\delta_j}=d^{(\frac{N}{2})}_{\delta_j}\le\frac{(\ell n)!}{2^\frac{\ell n}{2}(\frac{\ell n}{2})!}T\kappa^\ell\delta_j^{\alpha\ell}\le T\kappa^\ell\delta_j^{\alpha\ell}\sqrt{(\ell n)!}
	\]
by \eqref{dsd_comp_est} and \eqref{fact_est} which, together with \eqref{eq:rho}, yields that
	\begin{equation}
	\label{eq:est_1}
		Y_{\ell,\delta_j} \lesssim_{n,\kappa}(\bm{C}_q+T+1) \, \delta_j^\alpha(q-1)^\frac{n}{2}\ell^\frac{n}{2}
	\end{equation}
holds for every $j\ge j_0$, every even $\ell \geq 2$, and every $q>2$ almost surely. By interpolation,
	\begin{equation}
	\label{eq:norm_of_G}
		\|G(\cdot+\delta_j)-G(\cdot)\|_{L^p(0,T-\delta_j)}\lesssim_{n,\kappa}(\bm{C}_q+T+1)\, \delta_j^\alpha (q-1)^\frac{n}{2} p^\frac{n}{2}
	\end{equation}
is obtained for every $j\ge j_0$, every $p\ge 1$, and every $q>2$ almost surely and, consequently, 
	\[ [G]_{B^\alpha_{\varPhi_{2/n},\infty}(0,T)} \lesssim_{n,\kappa} (\bm{C}_q+T+1)(q-1)^\frac{n}{2}\]
holds for every $q>2$ almost surely and
	\[ \left[\E[G]_{B^\alpha_{\varPhi_{2/n},\infty}(0,T)}^q\right]^\frac{1}{q}\lesssim (q-1)^\frac{n}{2}\]
holds for every $q>2$. The first assertion then follows from the above estimate by noticing that $G(0)\in L^{\varPhi_{2/n}}(\Omega)$ holds by \cite[Corollary 2.8.14]{NouPec12} and appealing to \cite[Corollary 26]{Sim90}. If, additionally, process $G$ satisfies condition \ref{ass:G2}, it follows by Jensen's inequality and Fatou's lemma that 
	\[\liminf_{j\rightarrow\infty} \delta_j^{-\alpha\ell} \mathbb{E} Y_{\ell,\delta_j}^\ell \geq T{\kappa^\prime}^\ell(n!)^\frac{\ell}{2} >0\]
holds for every even $\ell \geq 2$ from which it follows, by appealing to \eqref{eq:converge_2}, that also
	\begin{equation*}
	\liminf_{j\rightarrow\infty} \delta_j^{-\alpha \ell} Y_{\ell,\delta_j}^\ell \geq T{\kappa^\prime}^\ell(n!)^\frac{\ell}{2}>0
	\end{equation*}
holds for every even $\ell \geq 2$ almost surely which yields that
	\begin{equation}
	\label{eq:est_2}
	\liminf_{j\rightarrow\infty} \delta_j^{-\alpha}Y_{\ell,\delta_j} \geq T^\frac{1}{\ell}\kappa^\prime \sqrt{n!}>0
	\end{equation}
holds for every even $\ell \geq 2$ almost surely. By interpolation, the inequality
	\[Y_{2,\delta_j} \lesssim Y_{1,\delta_j}^\frac{1}{3} Y_{4,\delta_j}^\frac{2}{3}\]
holds almost surely and, by using estimates \eqref{eq:est_1} and \eqref{eq:est_2}, we obtain that
	\[ \liminf_{j\rightarrow\infty} \delta_j^{-\alpha} Y_{1,\delta_j} \geom 4^{-n}(n!)^\frac{3}{2}>0\]
holds almost surely. As a consequence, we have for $q\geq 1$ that 
	\[ [G]_{B_{1,q}^\alpha(0,T)} = \left(\sum_{j\geq j_0} \delta_j^{-\alpha q}\|G(\cdot+\delta_j)-G(\cdot)\|_{L^1(0,T-\delta_j)}^q\right)^\frac{1}{q}\]
diverges almost surely so that the paths of $G$ do not belong to the space $B_{1,q}^\alpha(0,T)$. The second assertion of the theorem follows by the embedding of Besov spaces from \cite[Theorem 3.3.1]{Trie83}. \qed

\subsection{Proof of \autoref{main_thm} with condition \texorpdfstring{\ref{ass:G3'}}{(G3')}}
\label{sec:app_A}
The proof follows the same strategy as the proof in \autoref{sec:proof}. The difference comes in the estimate of $\| d_{\delta}^{(k)}\|_{H^{\otimes_2^{N-2k}}}$ for $k<N/2$. In particular, let $J_\delta\in\mathbb{N}$ be such that $\lambda_\delta = \frac{T-\delta}{J_\delta}\geq \delta$. Then we have 
\begin{align*}
		\MoveEqLeft[1] \|d^{(k)}_{\delta}\|_{H^{\otimes_2^{N-2k}}} &  \\
		& \leq\sum_{V\in\mathcal E^k_{I_1,\dots,I_\ell}}\left\|\int_0^{T-\delta}R^V(A_{s,\delta},\dots,A_{s,\delta})\d{s}\right\|_{H^{\otimes_2^{N-2k}}}
		\\
		&=\sum_{V\in\mathcal E^k_{I_1,\dots,I_\ell}}\left\|\lambda_\delta\sum_{j=1}^{J_\delta}\int_0^1R^V(A_{(j-1)\lambda_\delta+s\lambda_\delta,\delta},\dots,A_{(j-1)\lambda_\delta+s\lambda_\delta,\delta})\d{s}\right\|_{H^{\otimes_2^{N-2k}}}
		\\
		&\leq \sum_{V\in\mathcal E^k_{I_1,\dots,I_\ell}}\lambda_\delta\int_0^1\left\|\sum_{j=1}^{J_\delta}R^V(A_{(j-1)\lambda_\delta+s\lambda_\delta,\delta},\dots,A_{(j-1)\lambda_\delta+s\lambda_\delta,\delta})\right\|_{H^{\otimes_2^{N-2k}}}\d{s}
		\\
		&\leq\sum_{V\in\mathcal E^k_{I_1,\dots,I_\ell}}\lambda_\delta\int_0^1\Bigg(\sum_{m=1}^{J_\delta}\sum_{m^\prime=1}^{J_\delta}|\langle R^V(A_{(m-1)\lambda_\delta+s\lambda_\delta,\delta},\dots,A_{(m-1)\lambda_\delta+s\lambda_\delta,\delta}),
		\\ 
		& \hspace{5cm} R^V(A_{(m^\prime-1)\lambda_\delta+s\lambda_\delta,\delta},\dots,A_{(m^\prime-1)\lambda_\delta+s\lambda_\delta,\delta})\rangle_{H^{\otimes_2^{N-2k}}}|\Bigg)^{\frac12}\d{s}
		\\
		&\leq\sum_{V\in\mathcal E^k_{I_1,\dots,I_\ell}}\lambda_\delta\int_0^1\left(\sum_{m=1}^{J_\delta}\sum_{m^\prime=1}^{J_\delta}\prod_{i=1}^\ell\|R^{V_i}(A_{(m-1)\lambda_\delta+s\lambda_\delta,\delta},A_{(m^\prime-1)\lambda_\delta+s\lambda_\delta,\delta})\|_{H^{\otimes_2^{2n-2|V_i|}}}\right)^{\frac12}\d{s}
	\end{align*}
by using \autoref{mult_est_2}. Now, we have in fact that $k\leq N-2$ as $N$ is even which means that the set $V_*$ defined by \eqref{vstar} contains at least two elements. It follows that either there are $i,i^\prime\in\{1,\ldots,\ell\}$, $i\neq i^\prime$, such that $V_i\neq \emptyset$ and $V_{i^\prime}\neq \emptyset$ or there is $i\in\{1,\ldots,\ell\}$ such that $V_i$ contains at least two elements. Consequently, the estimate 
	\[ \prod_{i=1}^\ell \|R^{V_i}(A_{(m-1)\lambda_\delta + s\lambda_\delta,\delta},A_{(m^\prime-1)\lambda_\delta+s\lambda_\delta,\delta})\|_{H^{\otimes_2^{2n-2|V_i|}}} \leq \kappa^{2\ell} \delta^{2\alpha\ell} K(\delta,|m-m^\prime|\lambda_\delta)\]
for $m,m^\prime\in\{1,\ldots, J_\delta\}$, $m\neq m^\prime$, is obtained from assumption \ref{ass:G1} and the assumption that
	\[ C_{s,s} (x,y)\leq K(s, |x-y|)\]
holds for $x,y\geq 0$ and $s>0$ such that $x+s\leq T$, $y+s\leq T$, and $(x,x+s)\cap(y,y+s)=\emptyset$ by \autoref{mult_est}. This, together with \eqref{card_e}, yields
	\[ \|d_{\delta}^{(k)}\|_{H^{\otimes_2^{N-2k}}}\leq \frac{N!\kappa^{\ell}}{2^kk! (N-2k)!} \delta^{\alpha\ell} [\tilde{F}(\delta)]^\frac{1}{2}\]
where 
	\[ \tilde{F}(\delta) = \lambda_\delta^2 \left[ n\kappa^4J_\delta + \sum_{\substack{m, m^\prime=1 \\ m\neq m^\prime}}^{J_\delta} K(\delta,|m-m^\prime|\lambda_\delta)\right]\]
and hence
	\[ \mathbb{E}(Y_{\ell,\delta}^\ell - \mathbb{E} Y_{\ell,\delta}^\ell)^2 \leq \kappa^{2\ell}\delta^{2\alpha\ell} (\ell n)! 2^{\ell n} \tilde{F}(\delta).\]
as in \eqref{eq:second_moment}. It now follows for $q>2$ by \eqref{eq:hypercontractivity} that
	\[ 
		\E|Y_{\ell,\delta}^\delta-\E Y_{\ell,\delta}^\ell|^q 
			\leq
		 (\ell n)^\frac{q}{2}(q-1)^{\frac{q\ell n}{2}}\kappa^{\ell q}\delta^{\alpha\ell q}[(\ell n)!]^\frac{q}{2} 2^{\frac{\ell nq}{2}}[\tilde{F}(\delta)]^\frac{q}{2},\]
and if we define $\delta_j=T2^{-j}$, $J_{\delta_j}=2^j-1$ (so that $\lambda_{\delta_j}=\delta_j$), and 
	\[ \bm{C}_q=\left[ \sum_{j\in\mathbb{N}}\sum_{\ell\in 2\mathbb{N}}\left(\frac{|Y_{\ell,\delta_j}^\ell-\mathbb{E}Y_{\ell,\delta_j}^\ell|}{(2\kappa)^\ell (q-1)^\frac{\ell n}{2} \delta_j^{\alpha\ell} \sqrt{2^{\ell n}(\ell n)(\ell n)!}}\right)^q\right]^\frac{1}{q},\]
then  $[\mathbb{E}\bm{C}_q^q]^\frac{1}{q}\leq [\sum_{j\in\mathbb{N}} \tilde{F}(\delta_j)]^\frac{1}{2}$ where the finiteness of the sum follows from the assumption 
	\[ \sum_{j=1}^\infty \delta_j^{2}\sum_{\substack{m,m^\prime=1\\ m\neq m^\prime}}^{J_{\delta_j}} K(\delta_j, |m-m^\prime|\delta_j)<\infty.\]
The rest of the proof follows as in \autoref{sec:proof} (with $j_0=1$). \qed

\subsection{Proof of \autoref{thm:main_thm_2}}
\label{sec:proof_2}
We shall use the same notation as in the proof of \autoref{main_thm} in \autoref{sec:proof}. In there, the upper bounds are already obtained under assumption \ref{ass:G3}. (Indeed, the map $\delta\mapsto \|G(\cdot+\delta)-G(\cdot)\|_{L^p(0,T-\delta)}$ is sub-additive and lower semi-continuous so that it follows from \eqref{eq:norm_of_G} that $\|G(\cdot+\delta)-G(\cdot)\|_{L^p(0,T-\delta)}\leom \delta^\alpha p^\frac{n}{2}$ holds for every $\delta\in (0,T)$ and $p\geq 1$ almost surely.)  To obtain the lower bounds, it suffices to consider the case $p=1$. The lower bound for this case follows from \eqref{eq:est_2} by interpolation once it is shown that the process
 \[f_G^{\ell,\alpha}(r)=r^{-\alpha\ell}[Y_{\ell,r}^\ell-\mathbb{E} Y_{\ell,r}^\ell], \qquad r\in [0,T),\]
where $\ell\geq 2$ is a positive even integer, has a continuous version. To this end, let $\ell\geq 2$ be a positive even integer, define $N=\ell n$, and let $k$ be an integer such that $0\leq k<\frac{N}{2}$. Let also $0<s<t<T$. By proceeding similarly as in the proof of \autoref{main_thm}, we obtain 
	\begin{equation}
	\label{eq:mt_1}
	 \E [f_{G}^{\ell,\alpha}(t) - f_G^{\ell, \alpha}(s)]^2 \leq \sum_{k=0}^{\frac{N}{2}-1} (N-2k)! \|t^{-\alpha\ell}d_t^{(k)} - s^{-\alpha\ell}d_s^{(k)}\|_{H^{\otimes_2^{N-2k}}}^2.
	 \end{equation}
The norm of the difference can be estimated as
	\begin{align*}
		 \|t^{-\alpha\ell}d_t^{(k)} - s^{-\alpha\ell}d_s^{(k)}\|_{H^{\otimes_2^{N-2k}}} & \\
		 &\hspace{-4cm} \leq \sum_{V\in\mathcal{E}^k_{I_1, \ldots, I_{\ell}}} \left\|\int_0^{T-t}t^{-\alpha\ell} R^V(A_{x,t},\ldots,A_{x,t})\d{x} - \int_0^{T-s} s^{-\alpha\ell}R^V(A_{x,s},\ldots, A_{x,s})\d{x}\right\|_{H^{\otimes_2^{N-2k}}}
\intertext{and, upon denoting \[\bm{A}_{z,r}=(\underbrace{A_{z,r},\ldots,A_{z,r}}_{\ell\times})\] for $r\in (0,T)$ and $z\in (0,T-r)$ for simplicity, the chain continues as}
		& \hspace{-4cm} \leq  \sum_{V\in\mathcal{E}^k_{I_1, \ldots, I_{\ell}}} \Bigg\{\left\|\int_0^{T-t}[t^{-\alpha\ell} R^V(\bm{A}_{x,t}) - s^{-\alpha\ell}R^V(\bm{A}_{x,s})]\d{x}\right\|_{H^{\otimes_2^{N-2k}}} \\
		& \hspace{2cm} + \left\| \int_{T-t}^{T-s}s^{-\alpha\ell} R^V(\bm{A}_{x,s})\d{x}\right\|_{H^{\otimes_2^{N-2k}}}\Bigg\}.\numberthis \label{eq:mt_2}
	\end{align*}
Denote the first and the second term in the sum by $I_1$ and $I_2$, respectively. We have 
	\[ I_2 \leq s^{-\alpha\ell}(t-s)\sup_{x\in (T-t,T-s)}\|\bm{A}_{x,s}\|_{H^{\otimes_2^{N-2k}}} \leq s^{-\alpha\ell}(t-s)\sup_{x\in (T-t,T-s)} \|A_{x,s}\|_{H^{\otimes_2^n}}^{\ell} \leq \kappa^\ell (t-s)\]
by using \autoref{mult_est} and assumption \ref{ass:G1} successively. The focus is on term $I_1$ now. Write
	\begin{align*}
		I_1^2 & = t^{-\alpha\ell} \int_0^{T-t}\int_0^{T-t} \langle R^V(\bm{A}_{x,t}),R^V(\bm{A}_{y,t})\rangle_{H^{\otimes_2^{N-2k}}}\d{x}\d{y} \\
			&\hspace{1cm} -2(st)^{-\alpha\ell} \int_0^{T-t}\int_0^{T-t} \langle R^V(\bm{A}_{x,t}),R^V(\bm{A}_{y,s})\rangle_{H^{\otimes_2^{N-2k}}}\d{x}\d{y}\\
			&\hspace{1cm} + s^{-2\alpha\ell} \int_0^{T-t}\int_0^{T-t} \langle R^V(\bm{A}_{x,s}),R^V(\bm{A}_{y,s})\rangle_{H^{\otimes_2^{N-2k}}}\d{x}\d{y} \numberthis\label{eq:split}
	\end{align*}
Denote $\lambda=\frac{s}{t}$ for simplicity and assume for now that $\frac{t}{2}\leq s<t$. From \eqref{eq:split}, we also have 
	\begin{align*}
		I_1^2 & = t^{-2\alpha\ell}(1-\lambda^{-\alpha\ell})\int_0^{T-t}\int_0^{T-t} \langle R^V(\bm{A}_{x,t}),R^V(\bm{A}_{y,t})\rangle_{H^{\otimes_2^{N-2k}}}\d{x}\d{y} 
		\\
			& \hspace{1cm} + t^{-2\alpha\ell} (\lambda^{-2\alpha\ell} - \lambda^{-\alpha\ell})\int_0^{T-t}\int_0^{T-t} \langle R^V(\bm{A}_{x,s}),R^V(\bm{A}_{y,s})\rangle_{H^{\otimes_2^{N-2k}}}\d{x}\d{y} 
		\\
			& \hspace{1cm} + t^{-2\alpha\ell}\lambda^{-\alpha\ell} \int_0^{T-t}\int_0^{T-t} \langle R^V(\bm{A}_{x,t}) - R^V(\bm{A}_{x,s}),R^V(\bm{A}_{y,t})-R^V(\bm{A}_{y,s})\rangle_{H^{\otimes_2^{N-2k}}}\d{x}\d{y}.
	\end{align*}
Now, denote
	\[ \bm{A}_{z, r_1,r_2}^i=(A_{z,r_1},\ldots, A_{z,r_1},\underbrace{A_{z,r_1}-A_{z,r_2}}_{i^{\mbox{\tiny th}}\mbox{ \scriptsize position}},A_{z,r_2},\ldots,A_{z,r_2})\]
for $r_1,r_2\in (0,T),$ $z\in (0,T-\max\{r_1,r_2\})$, $i\in \{1,\ldots, \ell\}$; and denote also its $m$\textsuperscript{th} element, $m\in\{1,\ldots,\ell\}$, by $[\bm{A}_{z,r_1,r_2}^i]_m$. Then for $x,y\in (0,T-t)$, we have
\begin{align*}\label{i2est}
\left|\langle R^V(\bm{A}_{x,t})-R^V(\bm{A}_{x,s}), R^V(\bm{A}_{y,t})-R^V(\bm{A}_{y,s})\rangle_{H^{\otimes_2^{N-2k}}}\right| & \\
& \hspace{-8cm} \leq\sum_{i=1}^\ell\sum_{j=1}^\ell\left|\langle R^V(\bm{A}_{x,t,s}^i),R^V(\bm{A}_{y,t,s}^j)\rangle_{H^{\otimes_2^{N-2k}}}\right|\\
&\hspace{-8cm} \leq\sum_{i=1}^\ell\sum_{j=1}^\ell\left\{\|[\bm{A}_{x,t,s}^i]_{\tilde{m}}\otimes_l[\bm{A}_{y,t,s}^j]_{\tilde{m}}\|_{H^{\otimes_2^{2n-2l}}}\prod_{m\neq\tilde m}\left\|[\bm{A}_{x,t,s}^i]_{m}\right\|_{H^{\otimes_2^{n}}}\|[\bm{A}_{y,t,s}^j]_m\|_{H^{\otimes_2^{n}}}\right\}
\end{align*} 
by \autoref{mult_est} and \autoref{mult_est_2} since $V_{\tilde{m}}\neq\emptyset$ for at least one index $\tilde m\in\{1,\dots,\ell\}$, i.e. $l\in\{1,\dots,\ell\}$. Now, denote the product in the brackets above by $J^{i,j}_{t,s,\tilde{m}}(x,y) $ and define 
\[ \tilde{C}_{r_1,r_2}(z_1,z_2) = \sum_{j=1}^n \|A_{z_1,r_1}\otimes_jA_{z_2,r_2}\|_{H^{\otimes_2^{2(n-j)}}}\]
for $z_1,z_2\geq 0$ and $r_1,r_2>0$ that satisfy $z_1+r_1\leq T$, $z_2+r_2\leq T$. If $\tilde{m}$ is, for example, such that $[\bm{A}_{x,t,s}^i]_{\tilde{m}}=A_{x,t}$ and $[\bm{A}_{y,t,s}^j]_{\tilde{m}} = A_{y,t}$, we have the estimate
	\begin{align*}
		J^{i,j}_{t,s,\tilde{m}}(x,y) & \leq \tilde{C}_{t,t}(x,y)\|A_{x,t}\|^{i-2}_{H^{\otimes_2^n}}\|A_{x,s}\|^{\ell-i}_{H^{\otimes_2^n}}\|A_{x,t}-A_{x,s}\|_{H^{\otimes_2^n}}\\
		& \hspace{2cm} \cdot\|A_{y,t}\|^{j-2}_{H^{\otimes_2^n}}\|A_{y,s}\|^{\ell-j}_{H^{\otimes_2^n}} \|A_{y,t}-A_{y,s}\|_{H^{\otimes_2^n}}\\
		& \leq (\kappa^2t^{2\alpha})^{\ell-1}(1-\lambda)^{2\alpha}\tilde{C}_{t,t}(x,y)
	\end{align*}
by realizing that $A_{x,t}-A_{x,s} = A_{x+s,t-s}$ holds and using \ref{ass:G1}. By similar arguments, the estimates
	\begin{align*}
		J^{i,j}_{t,s,\tilde{m}}(x,y) & \leq (\kappa^2t^{2\alpha})^{\ell-1} 
		\begin{cases}
			(1-\lambda)^{2\alpha}\tilde C_{t,t}(x,y)\\
			(1-\lambda)^{2\alpha}\tilde C_{t,s}(x,y)\\
			(1-\lambda)^{2\alpha}\tilde C_{s,t}(x,y)\\
			(1-\lambda)^{2\alpha}\lambda^{-2\alpha}\tilde C_{s,s}(x,y)\\
			(1-\lambda)^{\alpha}\lambda^\alpha \tilde C_{t-s,t}(x+s,y)\\
			(1-\lambda)^{\alpha}\tilde C_{t-s,s}(x+s,y)\\
			(1-\lambda)^{\alpha}\lambda^\alpha \tilde C_{t,t-s}(x,y+s)\\
			(1-\lambda)^{\alpha}\tilde C_{s,t-s}(x,y+s)\\
			\tilde C_{t-s,t-s}(x+s,y+s)
		\end{cases}
	\end{align*}
depending on the precise value of $\tilde{m}$ can be obtained but in any case, we have the estimate 
	\[ \int_0^{T-t}\int_0^{T-t} J^{i,j}_{t,s,\tilde{m}}(x,y)\d{x}\d{y} \lesssim t^{2\alpha\ell + \varepsilon}(1-\lambda)^{2\alpha} \]
because assumption \ref{ass:G4} implies 
	\[ \int_0^{T-r_2}\int_0^{T-r_1} \tilde{C}_{r_1,r_2}(z_1,z_2)\d{z_1}\d{z_2} \lesssim  r_1^{\alpha+\frac{\varepsilon}{2}}r_2^{\alpha+\frac{\varepsilon}{2}}, \qquad r_1,r_2\in (0,T).\]
Thus, the inequality
	\[
		\int_0^{T-t}\int_0^{T-t} \left|\langle R^V(\bm{A}_{x,t}) - R^V(\bm{A}_{x,s}),R^V(\bm{A}_{y,t})-R^V(\bm{A}_{y,s})\rangle_{H^{\otimes_2^{N-2k}}}\right|\d{x}\d{y} \lesssim t^{2\alpha\ell + \varepsilon}(1-\lambda)^{2\alpha}
	\]
is shown. Similarly, the estimates 
	\begin{align*}
	\int_0^{T-t}\int_0^{T-t} \left| \langle R^V(\bm{A}_{x,s}),R^V(\bm{A}_{y,s})\rangle_{H^{\otimes_2^{N-2k}}}\right|\d{x}\d{y} & \lesssim s^{2\alpha\ell + \varepsilon}\\
		\int_0^{T-t}\int_0^{T-t} \left|\langle R^V(\bm{A}_{x,t}),R^V(\bm{A}_{y,t})\rangle_{H^{\otimes_2^{N-2k}}}\right|\d{x}\d{y} &  \lesssim t^{2\alpha\ell + \varepsilon}
	\end{align*}
are shown to hold by appealing to \autoref{mult_est}, \autoref{mult_est_2}, and assumptions \ref{ass:G1} and \ref{ass:G4}. Consequently, 
	\[I_1^2 \lesssim t^{\varepsilon}(\lambda^{-\alpha\ell}-\lambda^{\alpha\ell})+ t^{\varepsilon}\lambda^{-\alpha\ell}(1-\lambda)^{2\alpha} \lesssim t^{\varepsilon}(1-\lambda)^{\min\,\{1,2\alpha\}}\]
as $\frac{t}{2}\leq s<t$. Now, if $0<s<\frac{t}{2}$, proceeding as above via \autoref{mult_est} and \autoref{mult_est_2}, equality \eqref{eq:split} yields 
	\begin{align*}
		I_1^2 & \leq t^{-2\alpha\ell}(\kappa^2t^{2\alpha})^{\ell-1}\int_0^{T-t}\int_0^{T-t}\tilde{C}_{t,t}(x,y)\d{x}\d{y}\\
&\hspace{1cm} +2t^{-\alpha\ell}s^{-\alpha\ell}(\kappa^2 t^\alpha s^\alpha)^{\ell-1}\int_0^{T-t}\int_0^{T-t}\tilde{C}_{t,s}(x,y)\d{x}\d{y}\\
&\hspace{1cm} +s^{-2\alpha\ell}(\kappa^2s^{2\alpha})^{\ell-1}\int_0^{T-t}\int_0^{T-t}\tilde{C}_{s,s}(x,y)\d{x}\d{y}
	\end{align*}
so that 
	\[I_1^2 \lesssim t^{\varepsilon} \lesssim  t^{\varepsilon}(1-\lambda)\] 
holds by \ref{ass:G1} and \ref{ass:G4}. Thus we obtain 
	\[\|t^{-\alpha\ell}d_t^{(k)} - s^{-\alpha\ell}d_s^{(k)}\|_{H^{\otimes_2^{N-2k}}} \lesssim t^{\frac{\varepsilon}{2}}(1-\lambda)^{\min\{\frac{1}{2},\alpha\}}\]
from \eqref{eq:mt_2}. Consequently, it follows from \eqref{eq:mt_1} by Kolmogorov's continuity theorem that process $f_G^{\ell,\alpha}$ has a continuous version.\qed

\section*{Declarations}


\textbf{Funding information.} This research was supported by the Czech Science Foundation project No.~19-07140S. 

\medskip

\textbf{Conflict of interest.} The authors have no competing interests to declare that are relevant to the content of this article.

\medskip

\textbf{Data sharing.} Data sharing is not applicable to this article as no datasets were generated or analysed during the current work.


\end{document}